\documentclass[11pt]{article}
\usepackage{fullpage,graphicx,amsmath,amsthm,amssymb,lineno}
\def\bfm#1{\mbox{\boldmath$#1$}}
\def\qed{\hfill \rule{4pt}{7pt}}
\newtheorem{theorem}{Theorem}[section]
\newtheorem{lemma}[theorem]{Lemma}
\newtheorem{conjecture}[theorem]{Conjecture}
\newtheorem{corollary}[theorem]{Corollary}
\newtheorem{proposition}[theorem]{Proposition}
\numberwithin{equation}{section}
\numberwithin{figure}{section}
\begin{document}

\title{\bf On Box-Perfect Graphs}

\author{Guoli Ding\thanks{Supported in part by NSF grant DMS-1500699.} \\ Department of Mathematics, 
Louisiana State University, Baton Rouge, USA \and 
Wenan Zang\thanks{Supported in part by the Research Grants Council of Hong Kong.}  \,\, and \,\, 
Qiulan Zhao\thanks{Corresponding author. E-mail:  qiulanzhao@163.com.}\\
Department of Mathematics, The University of Hong Kong, Hong Kong, China}

\date{\today}

\maketitle

\begin{abstract}
Let $G=(V,E)$ be a graph and let $A_G$ be the clique-vertex incidence matrix of $G$. It is well known that
$G$ is perfect iff the system $A_{_G}\bfm x\le \bfm 1$, $\bfm x\ge\bfm0$ is totally dual integral (TDI).
In 1982, Cameron and Edmonds proposed to call $G$ box-perfect if the system $A_{_G}\bfm x\le \bfm 1$, 
$\bfm x\ge\bfm0$ is box-totally dual integral (box-TDI), and posed the problem of characterizing such 
graphs. In this paper we prove the Cameron-Edmonds conjecture on box-perfectness of parity graphs, and 
identify several other classes of box-perfect graphs. We also develop a general and powerful method for 
establishing box-perfectness. 

\end{abstract}

\section{Introduction}

A rational system $A\bfm x\le \bfm b$ is called {\it totally dual integral} (TDI) if the minimum in the 
LP-duality equation 
\begin{equation}\label{eq:lp}
\max\{\bfm w^T\bfm x: A\bfm x\le \bfm b\} = \min\{\bfm y^T\bfm b: \bfm y^TA= \bfm w^T; \ \bfm y\ge \bfm 0\}
\end{equation}
has an integral optimal solution, for every integral vector $\bfm w$ for which the minimum is finite. 
Edmonds and Giles \cite{tdi} proved that total dual integrality implies primal integrality: if $A\bfm x 
\le \bfm b$ is TDI and $\bfm b$ is integral, then both programs in (\ref{eq:lp}) have integral 
optimal solutions whenever they have finite optimum.  So the model of TDI systems serves as a general
framework for establishing min-max results in combinatorial optimization (see Schrijver \cite{schrijver3}
for an comprehensive and in-depth account).  As summarized by Schrijver \cite{S}, the importance of a 
min-max relation is twofold: first, it serves as an optimality criterion and as a good characterization 
for the corresponding optimization problem; second, a min-max relation frequently yields an elegant 
combinatorial theorem, and allows a geometrical representation of the corresponding problem in 
terms of a polyhedron.  Many well-known results and difficult conjectures in combinatorial optimization 
can be rephrased as saying that a certain linear system is TDI; in particular, by Lov\'asz' Replication 
Lemma \cite{L}, a graph $G$ is perfect if and only if the system $A_{_G}\bfm x\le \bfm 1$, $\bfm x\ge\bfm0$ 
is TDI, where $A_G$ is the clique-vertex incidence matrix of $G$. The reader is referred to 
Chudnovsky {\em et al.} \cite{CRST,ep} for the proof of the Strong Perfect Graph Theorem and
to Chudnovsky {\em et al.} \cite{CCLSV} for recognition of perfect graphs.

A rational system $A\bfm x\le \bfm b$ is called {\it box-totally dual integral} (box-TDI) 
if $A\bfm x\le \bfm b$, $\bfm l\le \bfm x\le \bfm u$ is TDI for all vectors $\bfm l$ and $\bfm u$, 
where each coordinate of $\bfm l$ and $\bfm u$ is either a  rational number or $\pm\infty$. By taking 
$\bfm l=-\bfm \infty$ and $\bfm u=\bfm \infty$ it follows that every box-TDI system must be TDI. Cameron 
and Edmonds \cite{cameron,cameron2} proposed to call a graph $G$ {\it box-perfect} if the system  
$A_{_G}\bfm x\le \bfm 1$, $\bfm x\ge\bfm0$ is box-TDI; they also posed the problem of characterizing such graphs.

We make some preparations before presenting an equivalent definition of box-perfect graphs. 
Let $G=(V,E)$ be a graph (all graphs considered in this paper are simple unless otherwise stated). For any $X\subseteq V$, let $G[X]$ denote the subgraph of $G$ induced by $X$. For any $v\in V$, let $N_G(v)$ denote the set of vertices incident with $v$. Members of $N_G(v)$ are called {\it neighbors} of $v$. By {\it duplicating} a vertex $v$ of $G$ we obtain a new graph $G'$ constructed as follows: we first add a new vertex $v'$ to $G$, which may or may not be adjacent to $v$, and then we join $v'$ to all vertices in $N_G(v)$.

As usual, let $\alpha(G)$ and $\chi(G)$ denote respectively the stable number and chromatic number of $G$. 
Let $\bar\chi(G)=\chi(\bar G)$, which is the clique cover number of $G$.
For any integer $q\ge1$, let  \medskip\\
\indent $\alpha_q(G)=\max\{|X|: X\subseteq V(G)$ with $\chi(G[X])\le q\}$, and \\ 
\indent $\bar\chi_q(G)=\min\{q\bar\chi(G-X)+|X|: X\subseteq V(G)\}$. \medskip\\ 
Notice that $\alpha_1=\alpha$ and $\bar\chi_1=\bar\chi$. A graph $G$ is called {\it $q$-perfect} if $\alpha_q(G[X])=\bar\chi_q(G[X])$ holds for all $X\subseteq V(G)$. This concept was introduced by Lov\'asz \cite{lovasz} as an extension of perfect graphs, since 1-perfect graphs are precisely perfect graphs. Let us call a graph {\it totally perfect} if it is $q$-perfect for all integers $q\ge1$. Lov\'asz pointed out that comparability graphs, incomparability graphs, and line graphs of bipartite graphs are totally perfect. However, $S_3$ is not 2-perfect, showing that a perfect graph does not have to be $q$-perfect when $q>1$.

\begin{figure}[ht]
%\centerline{\includegraphics[scale=0.58]{Figures/s3.eps}}
\centerline{\includegraphics[scale=0.58]{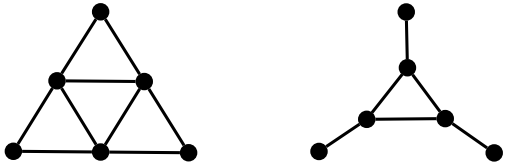}}
\caption{Graph $S_3$ and its complement $\bar S_3$}
\label{fig:s3}
\end{figure}

\begin{theorem}[Cameron \cite{cameron1}] \label{thm:cameron1} 
A graph is box-perfect if and only if every graph obtained from this graph by repeatedly duplicating vertices is totally perfect. 
\end{theorem}

This theorem implies the following immediately. 

\begin{corollary}[Cameron \cite{cameron1}] \label{cor:basic}
(1) Induced subgraphs of a box-perfect graph are box-perfect. \\ 
\indent (2) Duplicating vertices in a box-perfect graph results in a box-perfect graph. \\ 
\indent (3) Comparability and incomparability graphs are box-perfect. 
\end{corollary}

The next proposition contains a few other important observations made by Cameron  \cite{cameron1}. A matrix $A$ is {\it totally unimodular} if the determinant of every square submatrix of $A$ is 0 or $\pm1$. A $\{0,1\}$-matrix $A$ is {\it balanced} if none of its submtrices is the vertex-edge incidence matrix of an odd cycle. For each graph $G$, let $B_G$ be the submatrix of $A_G$ obtained by keeping only rows that correspond to maximal cliques of $G$. Let us call $G$ {\it totally unimodular} or {\it balanced} if $B_G$ is totally unimodular or balanced. It is worth pointing out that bipartite graphs and  their line graphs are totally unimodular, and every totally unimodular graph is balanced. In addition, as shown by Berge \cite{berge2}, all balanced graphs are totally perfect. 
Let $\bar S_3^+$ be obtained from the complement $\bar S_3$ of $S_3$ by adding a new vertex $v$ and joining $v$ to all six vertices of $\bar S_3$.

\begin{proposition}[Cameron \cite{cameron1}] \label{prop:Sn}
(1) $\bar S_3^+$ is not box-perfect. \\
\indent (2)  Totally unimodular graphs are box-perfect. \\  
\indent (3) Balanced graphs do not have to be box-perfect, shown by $\bar S_3^+$. \\
\indent (4) The complement of a box-perfect graph does not have to be box-perfect, shown by $\bar S_3$. \\
\indent (5) Box-perfectness is not preserved under taking clique sums, shown by $S_3$.
\end{proposition}

As we have seen, many nice properties of perfect graphs are not satisfied by box-perfect graphs. Another property of this kind is substitution: substituting a vertex of a box-perfect graph by a box-perfect graph does not have to yield a box-perfect graph, as shown by $\bar S_3^+$ (which is obtained by substituting a vertex of $K_2$ with $\bar S_3$).  To our knowledge, almost none of the 
known summing operations that preserve perfectness can carry over to box-perfectness -- this makes it extremely hard to obtain
a structural characterization of box-perfect graphs! 

At this point, the only known box-perfect graphs are totally unimodular graphs, comparability graphs, incomparability graphs, and $p$-comparability graphs (where $p\ge1$ and 1-comparability graphs are precisely comparability graphs) \cite{cameron, cameron2}. 
Cameron and Edmonds \cite{cameron} conjectured that every parity graph is box-perfect. In this paper we confirm this conjecture
and identify several other classes of box-perfect graphs, including claw-free box-perfect graphs.
In the next section we construct a class $\cal R$ of non-box-perfect graphs, from which we characterize box-perfect split graphs. It turns out that every minimal non-box-perfect graph that we know of is contained in a graph from $\cal R$. This observation 
raises the question: is it true that a graph $G$ is box-perfect if and only if $G$ does not contain any graph in $\cal R$ as an induced subgraph?  

In addition to structural description, the other difficulty with the study of box-perfect graphs lies in the lack of a
proper tool for establishing box-perfectness.  In section 3 we introduce a so-called ESP property, which is sufficient 
for a graph to be box-perfect. Although recognizing box-perfectness is an optimization problem, our approach based on 
the ESP property is of transparent combinatorial nature and hence is fairly easy to work with. For convenience,
we call a graph ESP if it has the aforementioned ESP property.  In the remainder of this paper, we shall establish several 
classes of box-perfect graphs by showing that they are actually ESP, including all classes obtained by Cameron 
\cite{cameron, cameron1, cameron2}.  We strongly believe that the ESP property is exactly the tool one needs for the study of box-perfect graphs. 

\begin{conjecture}
A perfect graph is box-perfect if and only if it is ESP if and only if it contains none of the members of $\cal R$ as an induced subgraph.   
\end{conjecture}

We close this section by mentioning a result on the complexity of recognizing box-perfect graphs. 

\begin{theorem}[Cook \cite{cook}]
The class of box-perfect graphs is in co-NP.
\end{theorem}

\section{A class of non-box-perfect graphs}

Let $S_n$ be the graph obtained from cycle $v_1v_2...v_{2n}v_1$ by adding edges $v_iv_j$ for all distinct even $i,j$. It was proved in \cite{cameron1} that $S_{2n+1}$ is not box-perfect for all $n\ge 1$. In this section we construct a class of non-box-perfect graphs, which include $\bar S_3^+$ and $S_{2n+1}$ ($n\ge 1$). We will use this result to characterize box-perfect split graphs (a graph is {\it split} if its vertex set can be partitioned into a clique and a stable set).

Let $G=(U,V,E)$ be a bipartite graph, where $U=\{u_1,...,u_m\}$ and $V=\{v_1,...,v_n\}$. The {\it biadjacency matrix} of $G$ is the $\{0,1\}$-matrix $M$ of dimension $m\times n$ such that $M_{i,j}=1$ if and only if $u_iv_j\in E$. Let $\cal Q$ be the set of bipartite graphs $G$ such that its biadjacency matrix $M$ is not totally unimodular but all submatrices of $M$ are. The following is a classical result of Camion.

\begin{lemma}[Camion \cite{camion}]\label{lem:camion}
Every graph $G=(U,V,E)$ in $\cal Q$ is Eulerian. In addition, $G$ satisfies $|U|=|V|$ and $|E|\equiv 2 \pmod 4$.
\end{lemma}

Let $\cal R$ be the class of graphs constructed as follows. 
Take a bipartite graph $G'=(U,V,E')\in \cal Q$ and a graph $G''=(V,E'')$ such that $N_{G'}(u)$ is a clique of $G''$ 
for all $u\in U$. 
Let $G=(U\cup V, E'\cup E'')$. 
If there exists $u\in U$ with $N_{G'}(u)=V$ then $G-u$ belongs to $\cal R$; otherwise $G$ belongs to $\cal R$.

\begin{figure}[ht]
%\centerline{\includegraphics[scale=0.55]{Figures/s3p.eps}}
\centerline{\includegraphics[scale=0.55]{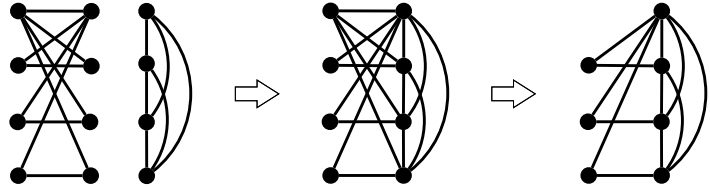}}
\caption{Graph $\bar S_3^+$ is constructed from a bipartite graph in $\cal Q$ and $K_4$}
\label{fig:s3p}
\end{figure}

\noindent{\bf Examples.} For each odd $n\ge3$, $S_n$ belongs to $\cal R$ since $S_n$ can be constructed from a cycle $G'=C_{2n}\in \cal Q$ and a complete graph $G''=K_n$, where no vertex is deleted in the construction. Graph $\bar S_3^+$ also belongs to $\cal R$. In this case a vertex is deleted in the construction, see Figure \ref{fig:s3p}.

\begin{lemma}\label{lem:O}
No graph in $\cal R$ is box-perfect.
\end{lemma}

\vspace{-2mm}
\noindent{\bf Proof.} Let $G\in\cal R$ be constructed from $G'=(U,V,E')\in \cal Q$ and $G''=(V,E'')$. Let $A_G$ and $B_G$ be the clique and maximal clique matrices of $G$. Then $A_G$ can be expressed as $A_G=[^{B_G}_{\ C}]$.
Let $M$ be the biadjacency matrix of $G'$ and let $n:=|U|\ (=|V|)$. Since every $u\in U$ belongs to exactly one maximal clique of $G$, the column of $B_G$ that corresponds to $u$ has precisely one nonzero entry. If no vertex was deleted in the construction of $G$ then $B_G$ can be expressed as $B_G=[^{M\ I_n}_{N\ \ \mathbf 0}]$, 
where the first $n$ columns are indexed by $V$ and the last $n$ columns are indexed by $U$. If a vertex $u_0\in U$ was deleted in the construction of $G$, then $G''$ has to be a complete graph. In this case, since $U$ does not have a second vertex adjacent to all vertices in $V$, $B_G$ can be expressed as $[M,J]$, where $J_{n\times (n-1)}=[^{I_{n-1}}_{\ \ \mathbf 0}]$ 
and the last row of $M$, which corresponds to $u_0$, is a vector of all ones. 
By Lemma \ref{lem:camion}, all entries of $\bfm1^TM$ and $M\bfm1$ are even, and $\bfm1^TM\bfm1=4m+2$, for an integer $m>0$. We consider the dual programs (with $A=A_G$)
\begin{equation}\label{eq:box}
\max\{\bfm w^T\bfm x: A\bfm x\le \bfm1; \bfm x\ge \bfm l\}=\min\{\bfm y^T\bfm1- \bfm z^T\bfm l: \bfm y^TA-\bfm z^T= \bfm w^T; \bfm y, \bfm z\ge\bfm0\}.
\end{equation}

Suppose no vertex was deleted in the construction of $G$. Let $p>2m+1$ be a prime and let 
$$\bfm w= \begin{bmatrix}\frac{1}{2}M^T\bfm 1\\ \bfm 0\end{bmatrix}, \ \ 
\bfm l=\begin{bmatrix}\bfm 0\\ \bfm1-\frac{1}{2p}M\bfm1\end{bmatrix}, \ \ 
\bfm x=\begin{bmatrix}\frac{1}{2p}\bfm1\\ \bfm1 - \frac{1}{2p} M \bfm1\end{bmatrix}, \ \ 
\bfm y=\begin{bmatrix}\frac{1}{2}\bfm 1\\ \bfm0\\ \bfm 0\end{bmatrix}, \ \ 
\bfm z=\begin{bmatrix}\bfm0\\ \frac{1}{2}\bfm 1\end{bmatrix}.$$ 
Then it is routine to verify that $\bfm w$ is integral, $\bfm l\ge\bfm0$, and $\bfm x, (\bfm y, \bfm z)$ are feasible solutions to (\ref{eq:box}). Moreover $\bfm w^T \bfm x =\frac{2m+1}{2p}= \bfm y^T\bfm1- \bfm z^T\bfm l$, so $\bfm x, (\bfm y, \bfm z)$ are optimal solutions. Since the optimal value is not $\frac{1}{p}$-integral, while $\bfm l$ is, it follows that the dual does not have an integral optimal solution and so $G$ is not box-perfect. Next, suppose that a vertex was deleted in the construction of $G$. The proof for this case is almost identical to the proof for the last case. The only difference is that $B_G$ has $2n-1$ columns, instead of $2n$ columns. Thus we need to truncate the corresponding vectors. To be precise, let 
$$\bfm w= \begin{bmatrix}\frac{1}{2}M^T\bfm 1\\ \bfm 0\end{bmatrix}, \ \ 
\bfm l=\begin{bmatrix}\bfm 0\\ J^T(\bfm1-\frac{1}{n}M\bfm1)\end{bmatrix}, \ \ 
\bfm x=\begin{bmatrix}\frac{1}{n}\bfm1\\ J^T(\bfm1 - \frac{1}{n} M \bfm1)\end{bmatrix}, \ \ 
\bfm y=\begin{bmatrix}\frac{1}{2}\bfm 1\\ \bfm0\end{bmatrix}, \ \ 
\bfm z=\begin{bmatrix}\bfm0\\ \frac{1}{2}\bfm 1\end{bmatrix}.$$
Using the fact that the last row of $M$ is $\bfm 1^T$ we deduce that $\bfm x$ and $(\bfm y, \bfm z)$ are feasible solutions, and $\bfm w^T \bfm x=\frac{2m+1}{n}= \bfm y^T\bfm 1- \bfm z^T \bfm l$, which implies that both solutions are optimal. Furthermore, since $M\bfm1$ is even and its last entry is $n$, we deduce that $n$ is even and thus $\bfm l$ is $\frac{1}{n/2}$-integral. However, the optimal value $\frac{2m+1}{n}$ is not $\frac{1}{n/2}$-integral, so $G$ is not box-perfect, which proves the theorem. \qed

To identify all minimally non-box-perfect split graphs, we consider the following subsets of $\cal Q$. Let $\mathcal Q_1$ consist of all bipartite graphs $G=(U,V,E)\in \cal Q$ such that $U$ has a vertex adjacent to all vertices of $V$. Let $\mathcal Q_2$ consist of all bipartite graphs $G=(U,V,E)\in \mathcal Q\backslash \mathcal Q_1$ such that the graph obtained from $G$ by adding a vertex and making it adjacent to all vertices of $V$ does not contain any graph in $\mathcal Q_1$ as an induced subgraph. Let $\cal S$ consist of all graphs in $\cal R$ that are constructed from a bipartite graph $G'\in \mathcal Q_1\cup \mathcal Q_2$ and a complete graph $G''$. It is clear that all members of $\cal S$ are split graphs. Moreover, $\bar S_3^+$ and $S_{2n+1}$ ($n\ge1$) belong to $\cal S$.

\begin{theorem} \label{thm:split}
The following are equivalent for any split graph $G$. \\ 
\indent (1) $G$ is box-perfect; \\ 
\indent (2) no graph in $\cal S$ is an induced subgraph of $G$; \\ 
\indent (3) $G$ is totally unimodular.
\end{theorem}

\vspace{-2mm}
\noindent{\bf Proof.} Implication (3) $\Rightarrow$ (1) follows from Proposition \ref{prop:Sn}(2) and implication (1) $\Rightarrow$ (2) follows from Lemma \ref{lem:O} and Corollary \ref{cor:basic}(1). To prove (2) $\Rightarrow$ (3), let $G=(U,V,E)$ be a split graph, where $U$ is a stable set and $V$ is a clique. Let $G''=G[V]$ and $G'=G\backslash E(G'')$. Let $G'''$ be the bipartite graph obtained from $G'$ by adding a vertex $w$ adjacent to all vertices in $V$. Let $M$ be the biadjacency matrix of $G'''$. 

We first prove that $M$ is totally unimodular. Suppose otherwise. Then $G'''$ has an induced subgraph $H'\in \cal Q$. Let us choose $H'$ so that $H'$ contains the new vertex $w$ whenever it is possible. Consequently, $H'\in\mathcal Q_1\cup\mathcal Q_2$. Let $H$ be constructed from $H'$ and a complete graph $H''$. Then $H\in\cal S$ and, by the construction of $G$, $G$ contains $H$ as an induced subgraph. This contradicts (2) and thus $M$ has to be totally unimodular.  

Let $N$ be the biadjacency matrix of $G'$. Then $B_G=[N,I]$ or $[^{N\ I}_{\,\mathbf 1\ \,\mathbf 0}]$, 
depending on if $V$ is a maximal clique of $G$. Notice that $M=[^N_{\, \mathbf 1}]$. So $B_G$, and thus $G$, is totally unimodular. \qed

This theorem shows that all minimally non-box-perfect split graphs are contained in $\cal S$. In fact, $\cal S$ consists of precisely such graphs.

\begin{theorem}\label{thm:minsplit}
A split graph $G$ belongs to $\cal S$ if and only if $G$ is not box-perfect but all its induced subgraphs are.
\end{theorem}

\vspace{-2mm}
\noindent{\bf Proof.} The backward implication follows immediately from Theorem \ref{thm:split}. To prove the forward implication, let $G\in\cal S$. By Lemma \ref{lem:O}, we only need to show that $G-w$ is box-perfect for all $w\in V(G)$. Suppose $G$ is constructed from a bipartite graph $G'=(U,V,E')\in \mathcal Q_1\cup \mathcal Q_2$ and a complete graph $G''=(V,E'')$. Let $M$ be the biadjacency matrix of $G'$ and let $n:=|U|=|V|$. Observe that if $G'\in\mathcal Q_1$ then $B_G=[N,J]$, where $N=M$ and $J=[^{I_{n-1}}_{\ \ \mathbf 0}]$; if $G'\in \mathcal Q_2$ then $B_G=[N,J]$, where $N=[^M_{\, \mathbf 1}]$ and $J=[^{I_n}_{\, \mathbf 0}]$.

Now it is straightforward to verify that, for each $u\in U$, $B_{G-u}=[N',J']$ is obtained from $B_G$ by deleting the row and the column indexed by $u$; for each $v\in V$, $B_{G-v}=[N',J']$ is obtained from $B_G$ by deleting the column indexed by $v$ and also possibly the last row. In both cases, $N'$ is a proper submatrix of $N$. This implies that $N'$ is totally unimodular and thus so is $[N',J']$. Consequently, $G-w$ is box-perfect (totally unimodular) for all $w\in V(G)$, which proves the theorem. \qed

\begin{figure}[ht]
%\centerline{\includegraphics[scale=0.5]{Figures/mingraph.eps}} 
\centerline{\includegraphics[scale=0.5]{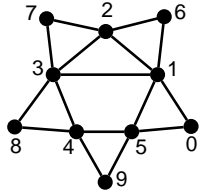}} 
\caption{A new non-box-perfect graph $G$}
\label{fig:newnbp}
\end{figure}

As we observed earlier that $\bar S_3^+$ and $S_{2n+1}$ $(n\ge 1)$ belong to $\cal S$. Thus these graphs are minimally non-box-perfect. We point out that, in addition to graphs in $\cal S$, other minimally non-box-perfect graphs can also be obtained using Lemma \ref{lem:O}. For instance, the graph illustrated in Figure \ref{fig:newnbp} is constructed from $G'=C_{10}$ and $G''=C_5+e$. By Lemma \ref{lem:O}, this graph $G$ is not box-perfect. However, $G$ is not minimally non-box-perfect since $H=G-\{9,0\}$ is not box-perfect, which is certified by vectors $\bfm w^T=(1, 1, 1, 1, 1, 0, 0, 0)$, $\bfm l^T=(0,0,0,0,0,\frac{1}{2},\frac{1}{2},\frac{1}{2})$, $\bfm x^T=(\frac{1}{4},\frac{1}{4},\frac{1}{4},\frac{1}{4},\frac{3}{4},\frac{1}{2},\frac{1}{2},\frac{1}{2})$, $\bfm y^T=(0,\frac{1}{2},\frac{1}{2},\frac{1}{2},\frac{1}{2},\frac{1}{2})$, and $\bfm z^T=(0,0,0,0,0,\frac{1}{2},\frac{1}{2},\frac{1}{2})$, where the first row of $B_H$ is the triangle 123.
It can be shown that $H$ is in fact minimally non-box-perfect because $H-x$ is totally unimodular for $x=1,2,3,4,5,6,7$, and $H-8$ has the ESP property defined in the next section which implies the box-perfectness.

\section{ESP graphs}

In this section we introduce a so-called ESP property, which is sufficient for a graph to be box-perfect. We shall
use this combinatorial property to identify several new classes of box-perfect graphs. We begin with a few lemmas.

\begin{lemma}[Chen, Ding and Zang \cite{cdz}] \label{lem:cdz}
Suppose $\bfm a_1$ and $\bfm a_2$ are rational vectors with $\bfm a_1 \ge \bfm a_2$, and $b_1$ and $b_2$ are rational numbers with $b_1 \le b_2$. Then the system $A\bfm x\le \bfm b$, $\bfm a^T_1\bfm x\le b_1$, $\bfm a^T_2\bfm x \le b_2$, $\bfm x \ge \bfm 0$ is box-TDI if and only if the system $A\bfm x\le \bfm b$, $\bfm a^T_1\bfm x\le b_1$, $\bfm x \ge \bfm 0$ is box-TDI.
\end{lemma}

\begin{lemma}[Cameron \cite{cameron1}] \label{lem:upper}
The system $A\bfm x\le \bfm b$ is box-TDI if and only if the system $A\bfm x\le \bfm b$, $\bfm x\le \bfm u$ is TDI, for all vectors $\bfm u$, where each coordinate of $\bfm u$ is either a rational number or $+\infty$.
\end{lemma}

The next two lemmas are reformulations of Theorem 22.7 and Theorem 22.13 of Schrijver \cite{schrijver}. 

\begin{lemma}[Schrijver \cite{schrijver}] \label{lem:inf}
Suppose the system $A\bfm x\le \bfm b$, $x_1\le u$ is TDI for all rational numbers $u$, where $x_1$ is the first coordinate of $\bfm x$. Then $A\bfm x\le \bfm b$ is TDI.
\end{lemma}

\begin{lemma}[Schrijver \cite{schrijver}] \label{lem:int}
A rational system $A\bfm x\le \bfm b$, $\bfm x \ge \bfm 0$ is TDI if and only if $\min\{\bfm y^T\bfm b: \bfm y^TA\ge \bfm w^T$\!, $\bfm y\ge\bfm 0$ {\rm is half-integral}\} is finite and is attained by an integral $\bfm y$, for each integral vector $\bfm w$ for which $\min\{\bfm y^T\bfm b: \bfm y^TA\ge \bfm w^T$\!, $\bfm y\ge\bfm 0\}$ is finite.
\end{lemma}

The next are two easy corollaries.  

\begin{lemma} \label{lem:A2B}
A graph $G$ is box-perfect if and only if the system $B_{_G}\bfm x\le\bfm 1$, $\bfm 0\le \bfm x\le \bfm u$ is TDI for all rational vectors $\bfm u\ge\bfm0$.
\end{lemma}

\vspace{-2mm}
\noindent{\bf Proof.} The forward implication follows immediately from the definition of box-TDI and Lemma \ref{lem:cdz}. Conversely, Lemma \ref{lem:upper} and Lemma \ref{lem:inf} imply that $B_{_G}\bfm x\le\bfm 1$, $\bfm x\ge \bfm 0$ is box-TDI. Then the result follows from Lemma \ref{lem:cdz}. \qed

\begin{lemma} \label{lem:halfi}
A graph $G$ is box-perfect if and only if for all rational $\bfm u\ge\bfm 0$ and integral $\bfm w\ge\bfm 0$, \medskip\\ 
\makebox[30.5mm]{} $\min\{\bfm  y^T\bfm 1 + \bfm  z^T\bfm u|\ \bfm  y^TB_{_G} + \bfm  z^T \ge 2\bfm w^T; \bfm  y, \bfm  z \ge\bfm0$ {\rm integral}\} \\ 
\makebox[24mm]{} $\ge 2\min\{\bfm  y^T\bfm 1 + \bfm  z^T\bfm u|\ \bfm  y^TB_{_G} + \bfm  z^T \ge \bfm w^T; \bfm  y, \bfm  z \ge\bfm0$ {\rm integral}\}. \hfill $(3.1)$
\end{lemma}

\vspace{-2mm}
\noindent {\bf Proof.} Observe that, for all vectors $\bfm u\ge\bfm 0$ and $\bfm w$, the three programs \smallskip\\ 
\indent\indent $\min\{\bfm  y^T\bfm 1 + \bfm  z^T\bfm u|\ \bfm  y^TB_{_G} + \bfm  z^T \ge \bfm w^T; \bfm  y, \bfm  z \ge\bfm0\}$ \\ 
\indent\indent $\min\{\bfm  y^T\bfm 1 + \bfm  z^T\bfm u|\ \bfm  y^TB_{_G} + \bfm  z^T \ge \bfm w^T; \bfm  y, \bfm  z \ge\bfm0$ half-integral\}\\ 
\indent\indent $\min\{\bfm  y^T\bfm 1 + \bfm  z^T\bfm u|\ \bfm  y^TB_{_G} + \bfm  z^T \ge \bfm w^T; \bfm  y, \bfm  z \ge\bfm0$ integral\}\smallskip\\ 
are finite. Moreover, replacing $\bfm w$ by $\bfm w_+$ does not change the minimum values of these programs, where $\bfm w_+$ is obtained from $\bfm w$ by turning its negative coordinates into zero. Therefore, the result follows immediately from Lemma \ref{lem:A2B} and Lemma \ref{lem:int}. \qed

Let $G=(V,E)$ be a graph. For any multiset $\Lambda$ of cliques of $G$ and any $v\in V$, let $d_{\Lambda}(v)$ denote the number of members of $\Lambda$ that contain $v$. We call $G$ {\it equitably subpatitionable} ({\it ESP}) if for every set $\Lambda$ of maximal cliques of $G$ there exist two multisets $\Lambda_1$ and $\Lambda_2$ of cliques of $G$ (which are not necessarily members of $\Lambda$) such that \smallskip\\ 
\indent (i) \ $|\Lambda_1| + |\Lambda_2| \le |\Lambda|$; \\ 
\indent (ii) \ $d_{\Lambda_1}(v)+ d_{\Lambda_2}(v)\ge d_{\Lambda}(v)$, for all $v\in V$; and \\ 
\indent (iii) \ $\min\{d_{\Lambda_1}(v),d_{\Lambda_2}(v)\}\ge \lfloor d_{\Lambda}(v)/2\rfloor$, for all $v\in V$. \smallskip\\ We call $(\Lambda_1,\Lambda_2)$ an {\it equitable subpartition} of $\Lambda$, and refer to the above (i), (ii), and (iii)
as {\em ESP property}. Note that (i) is equivalent to $|\Lambda_1| + |\Lambda_2| = |\Lambda|$ since we may include empty cliques in $\Lambda_1$ and $\Lambda_2$. Similarly, (ii) is equivalent to $d_{\Lambda_1}(v)+ d_{\Lambda_2}(v)= d_{\Lambda}(v)$ for all $v$, since cliques in $\Lambda_1,\Lambda_2$ can be replaced by smaller ones. Finally, it is also easy to see that in an ESP graph every multiset $\Lambda$ of cliques admits an equitable subpartition. We will use these facts without further explanation. 

\begin{theorem}\label{thm:esp}
Every ESP graph $G=(V,E)$ is box-perfect.
\end{theorem}

\vspace{-2mm}
\noindent{\bf Proof.} By Lemma \ref{lem:halfi} we only need to show that inequality (3.1) holds for all rational $\bfm u\ge \bfm0$ and all integral $\bfm w\ge \bfm 0$. Let $(\bfm  y^T, \bfm  z^T)$ be an optimal solution of the first minimum in (3.1). Let $\cal C$ be the set of maximal cliques of $G$ and let $\cal D$ be the multiset of members of $\cal C$ such that each $C\in\cal C$ appears in $\cal D$ exactly $y_C$ times. Let $\Lambda$ be the set of $C\in\cal C$ such that $y_C$ is odd.  Since $G$ is ESP, $\Lambda$ admits a equitable subpartition $(\Lambda_1,\Lambda_2)$. Since every clique can be extended into a maximal clique, we may assume without loss of generality that members of $\Lambda_1$ and $\Lambda_2$ are all in $\cal C$. Let $\mathcal D_0$ be the multiset of members of $\cal C$ such that each $C\in\cal C$ appears $\lfloor  y_C/2\rfloor$ times. It follows that $\mathcal D=\mathcal D_0\uplus \mathcal D_0\uplus \Lambda$, where $\uplus$ stands for multiset sum. For $i=1,2$, let $\mathcal D_i=\mathcal D_0\uplus \Lambda_i$. 
We deduce from (i) that 

\noindent(1) \ $|\mathcal D_1| + |\mathcal D_2|\le |\mathcal D|$.

Let $\bfm p = \bfm y^TB_{_G} + \bfm z^T - 2\bfm w^T$ and let $v\in V$. Without loss of generality we assume 

\noindent (2) \ $d_{\mathcal D_1}(v)\ge d_{\mathcal D_2}(v)$ \ and \ $\bfm p_v \bfm z_v=0$.

Since $d_{\mathcal D}(v)=2d_{\mathcal D_0}(v)+d_{\Lambda}(v)$, we deduce from (ii-iii) that $d_{\mathcal D_1}(v)+ d_{\mathcal D_2}(v) \ge d_{\mathcal D}(v)$ and $d_{\mathcal D_i}(v) = d_{\mathcal D_0}(v)+d_{\Lambda_i}(v) \ge \lfloor d_{\mathcal D}(v)/2\rfloor$ ($i=1,2$). Thus we conclude from (2) that

\noindent (3) \ $d_{\mathcal D_1}(v) \ge \lceil d_{\mathcal D}(v)/2\rceil$ \ and \ $d_{\mathcal D_2}(v) \ge \lfloor d_{\mathcal D}(v)/2\rfloor$. 

By the definition of $\cal D$ we have $d_{\cal D}(v)=\bfm  y^T B_v$, where $B_v$ is the column of $B_G$ indexed by $v$. So 

\noindent(4) \ $d_{\cal D}(v)+ \bfm z_v = \bfm p_v + 2\bfm w_v \ge 2\bfm w_v$. 

Since $\bfm w_v$ is an integer, we deduce that 

\noindent(5) \ $\bfm w_v\le \lfloor (d_{\cal D}(v)+ \bfm z_v)/2\rfloor$.

Setting $\bfm z_{1v}= \lfloor \bfm z_v/2\rfloor$ and $\bfm z_{2v}= \lceil \bfm z_v/2\rceil$, we have 

\noindent(6) \ $\bfm z_{1v}\bfm u_v + \bfm z_{2v} \bfm u_v= \bfm z_v\bfm u_v$. 

We further claim that 

\noindent(7) \ $d_{\mathcal D_i}(v) + \bfm z_{iv}\ge \bfm w_v$, for $i=1,2$. 

To see (7), recall $\bfm p_v \bfm z_v=0$ from (2). 
If $d_{\cal D}(v)$ is even, we deduce from (4) that $\bfm z_v$ is even, which implies, by (3-4), that $d_{\mathcal D_i}(v) + \bfm z_{iv}\ge \frac{1}{2} (d_{\mathcal D}(v) + \bfm z_v) \ge \bfm w_v$. 
So we assume that $d_{\cal D}(v)$ is odd. 
If $\bfm z_v=0$ then, by (3) and (5), $d_{\mathcal D_i}(v) + \bfm z_{iv} = d_{\mathcal D_i}(v) \ge \lfloor d_{\mathcal D}(v)/2\rfloor \ge \bfm w_v$. 
Else, by (2) and (4), $\bfm z_v$ is odd. Thus $d_{\mathcal D_i}(v) + \bfm z_{iv} \ge \frac{1}{2}(d_{\mathcal D_i}(v)\pm1)+\frac{1}{2}(\bfm z_v\mp1) = \frac{1}{2}(d_{\mathcal D_i}(v) + \bfm z_v)\ge \bfm w_v$, because of (3), (5), and the definition of $\bfm z_{iv}$. So (7) holds.

For $i=1,2$, let $\bfm  z_i=(\bfm z_{iv}: v\in V)$ and $\bfm y_i\in\mathbb Z_+^{\cal C}$ be the multiplicity function of $\mathcal D_i$. It follows from (7) that $\bfm  y_i^TB_{_G} + \bfm  z_i^T \ge \bfm w^T$, which means that both $(\bfm y_1,\bfm z_1)$ and $(\bfm y_2,\bfm z_2)$ are feasible solutions of the second program in (3.1). From (1) and (6) we also conclude that $\bfm  y_i^T\bfm 1 + \bfm  z_i^T\bfm u \le (\bfm  y^T\bfm 1 + \bfm  z^T\bfm u)/2$ holds for at least one $i\in\{1,2\}$. Hence inequality (3.1) holds, which proves the Theorem. \qed

For a perfect graph $G$, being ESP can be characterized as follows. Let $\mathbb Z_+$ denote the set of nonnegative integers. For any $d\in \mathbb Z_+^{V(G)}$, let $G^d$ denote the graph obtained from $G$ by substituting each vertex $v$ with a stable set of size $d(v)$. Note that $v$ is deleted when $d(v)=0$. Let $c_{_G} =\bfm 1^TB_G$. In other words, for each $v\in V(G)$, $c_{_G}(v)$ is the number of maximal cliques of $G$ that contain $v$. 

\begin{theorem}\label{thm:esp1}
Let $G$ be perfect. Then $G$ is ESP if and only if for every $d\in \mathbb Z_+^{V(G)}$ with $d\le c_{_G}$ there exists $d'\in \mathbb Z_+^{V(G)}$ such that $\lfloor d/2\rfloor \le d'\le \lceil d/2\rceil$ and $\alpha(G^{d'}) + \alpha(G^{d-d'}) \le  \alpha(G^d)$.
\end{theorem}

\vspace{-2mm}
\noindent{\bf Proof.} To prove the forward implication, let $G$ be ESP and let $d\in \mathbb Z_+^{V(G)}$. Since $G^d$ is perfect, its vertex set can be partitioned into $\alpha(G^d)$ cliques. These cliques naturally correspond to a multiset $\Lambda$ of $\alpha(G^d)$ cliques of $G$. Note that $|\Lambda|=\alpha(G^d)$ and $d_\Lambda = d$. Since $G$ is ESP, $\Lambda$ admits a equitable subpartition $(\Lambda_1,\Lambda_2)$. By deleting vertices from cliques in $\Lambda_1$ and $\Lambda_2$ we can obtained multisets $\Lambda_1^*$ and $\Lambda_2^*$ of cliques of $G$ such that $|\Lambda_1^*| + |\Lambda_2^*|\le |\Lambda_1| + |\Lambda_2|$, $d_{\Lambda_1^*} + d_{\Lambda_2^*} =d$, and $\min\{d_{\Lambda_1^*}, d_{\Lambda_2^*}\}\ge \lfloor d/2\rfloor$. Let $d'=d_{\Lambda_1^*}$. Then $\lfloor d/2\rfloor \le d'\le \lceil d/2\rceil$ and \smallskip\\ 
\indent\indent $\alpha(G^{d'}) + \alpha(G^{d-d'}) \le  \alpha(G^{d_{\Lambda_1}}) + \alpha(G^{d_{\Lambda_2}}) \le |\Lambda_1| + |\Lambda_2| \le |\Lambda| = \alpha(G^d)$, \smallskip\\ 
which proves the forward implication.

To prove the backward implication, let $\Lambda$ be a set of maximal cliques of $G$. Then $d:=d_{\Lambda}\le c_{_G}$ and thus there exists $d'$ as stated in the theorem. Let $d_1=d'$ and $d_2=d-d'$. For $i=1,2$, vertices of $G^{d_i}$ can be partitioned into $\alpha(G^{d_i})$ cliques, and these cliques correspond to a multiset $\Lambda_i$ of $\alpha(G^{d_i})$ cliques of $G$. Note that $d_{\Lambda_i} =d_i$. Thus $(\Lambda_1,\Lambda_2)$ is an equitable subpartition of $\Lambda$, which proves the theorem. \qed

We first remark that $\alpha(G^d)$ is exactly the maximum of $\sum_{v\in S} d(v)$ over all stable sets $S$ of $G$. Sometimes this interpretation is more convenient. We also remark that we do not know a box-perfect graph that is not ESP. It seems reasonable to conjecture that no such a graph exists. 

\section{Known box-perfect graphs}

Cameron \cite{cameron1} identified a few classes of box-perfect graphs. In this section we prove that they are in fact ESP graphs. Our results could be stronger than the results of Cameron if ESP and box-perfect are not equivalent. But the main reason for establishing our results is for future applications. We envision that more ESP graphs (possibly all box-perfect graphs) can be constructed from basic ESP graphs. Therefore, it is important to make sure that all known box-perfect graphs are ESP. 

\subsection{Totally unimodular graphs}

It is well known (see Theorem 19.3 of \cite{schrijver}) that in a totally unimodular matrix, each set of rows can be partitioned so that the sum of one part minus the sum of the other part is a $\{0,\pm1\}$-vector. If $G$ is totally unimodular then $B_G$ has this partition property, which implies immediately that $G$ satisfies the definition of ESP graphs. Thus we have the following.

\begin{theorem}
Totally unimodular graphs are ESP.
\end{theorem}

We point out that totally unimodular graphs include graphs like interval graphs, bipartite graphs, and block graphs (every block is a complete graph).

\subsection{Incomparability graphs}

\begin{theorem}\label{thm:incomp}
Every incomparability graph $G$ is ESP.
\end{theorem}

\vspace{-2mm}
\noindent{\bf Proof.} Since $G$ is perfect, we may apply Theorem \ref{thm:esp1}. 
Let $d\in\mathbb Z_+^{V(G)}$. Note that $G^d$ is again an incomparability graph. In fact, let $P$ be a poset such that $G$ is the incomparability of $P$ and let $P^d$ be obtained from $P$ by replacing each element $v$ with a chain of size $d(v)$. Then $G^d$ is the incomparability graph of poset $P^d$. For each positive integer $i$, let $A_i$ be the set of maximal elements of $P^d - (A_1\cup ... \cup A_{i-1})$. Then $(A_1, ..., A_n)$ is a partition of $V(G^d)$ into cliques, where $n=\alpha(G^d)$. Let $V_1$ be the union of $A_i$ for all odd $i$ and let $V_2$ be the union of $A_i$ for all even $i$. Then $G^d[V_1]$ and   $G^d[V_2]$ can be expressed as  $G^{d_1}$ and  $G^{d_2}$, respectively, for some $d_1,d_2\in \mathbb Z_+^{V(G)}$. It is easy to see that $d_1+d_2=d$ and $\lfloor d/2\rfloor \le d_j\le \lceil d/2\rceil$ ($j=1,2$). Moreover, each $\alpha(G^{d_j})$ is bounded by the number of $A_i$s contained in $V_j$. Therefore, $\alpha(G^{d_1}) + \alpha(G^{d_2}) \le  \alpha(G^d)$, which implies that $d'=d_1$ satisfies Theorem \ref{thm:esp1} and thus $G$ is ESP. \qed

\subsection{$p$-Comparability graphs}

$p$-Comparability graphs were introduced in \cite{cameron} and were shown \cite{cameron, cameron2} to be box-perfect. We show that they are ESP. Let $D$ be a digraph with a special set $T$ of vertices such that every arc is in a dicycle (directed cycle) and every dicycle meets $T$ exactly once. In particular, $D$ has no arc between any two vertices of $T$. If $p$ is an integer with $|T|\le p$, then a {\it p-comparability graph} $G$ is defined from $D$ by adding all chords of all dicycles, then deleting $T$, and finally ignoring all directions on edges. Note that 1-comparability graphs are precisely comparability graphs.

\begin{theorem} \label{thm:pcom}
Every $p$-comparability graph $G$ is ESP.
\end{theorem}

To prove this theorem we will need the following Lemma. Let $D=(V,A)$ be a digraph. For each dicycle $C$ of $D$, the {\it incidence vector} of $C$ is the vector $\chi^C\in \{0,1\}^A$ such that $\chi^C(a)=1$ if and only if $a$ is on $C$. A sum of incidence vectors of (not necessarily distinct) dicycles of $D$ is called a {\it circulation} of $D$. The following is a special case of  Corollary 11.2b of \cite{schrijver3}.

\begin{lemma}\label{lem:cir}
Every circulation $f$ is the sum of two circulations $f_1$, $f_2$ such that $\lfloor f/2\rfloor \le f_i \le \lceil f/2\rceil$ holds for both $i=1,2$. 
\end{lemma}

\noindent{\bf Proof of Theorem \ref{thm:pcom}.} Let $G$ be constructed from $D$ and $T$. Let $D^*$ be obtained from $D$ by splitting each vertex $v$ into $v'$ and $v''$ such that arcs entering $v$ are now entering $v'$, and arcs leaving $v$ are now leaving $v''$. We also add an arc from $v'$ to $v''$. Observe that for every dicycle $C$ of $D$, $D^*$ has a unique dicycle $C^*$ such that $A(C^*)\cap A(D) = A(C)$. Moreover, every dicycle of $D^*$ can be expressed as $C^*$ for a dicycle $C$ of $D$. 

We will use a fact proved in \cite{cameron2} that for every clique $K$ of $G$, there exists a dicycle $C_K$ of $D$ such that $K\subseteq V(C_K)$.

Let $\Lambda$ be a set of maximal cliques of $G$. We prove the theorem by showing that $\Lambda$ admits an equitable subpartition. Let $f$ be the sum of incidence vectors of $C^*_K$ over all $K\in \Lambda$. Since each $C_K$ meets $T$ exactly once, each $C_K^*$ must meet $T^*=\{t't'':t\in T\}$ exactly once. As a result, $|\Lambda|$ equals the sum of $f(a)$ over all $a\in T^*$. In addition, since each $K\in\Lambda$ is a maximal clique, we must have $V(C_K)-T=K$. This implies that $d_{\Lambda}(v) = f(v'v'')$ holds for all $v\in V(D)$.

Let $f_1$ and $f_2$ be the two circulations of $D^*$ determined by Lemma \ref{lem:cir}. For $i=1,2$, let $\mathcal C_i^*$ be the multiset of dicycles of $D^*$ such that $f_i$ is the sum of $\chi^{_{C^*}}$ over all $C^*\in\mathcal C_i^*$. Then let $\mathcal C_i$ be the multiset $\{C: C^*\in \mathcal C_i^*\}$ and $\Lambda_i=\{V(C)-T:C\in\mathcal C_i\}$. By the construction of $G$, each member of $\Lambda_i$ is a clique of $G$. Moreover, $d_{\Lambda_i}(v)=f_i(v'v')$ holds for all $v\in V(G)$, and $|\Lambda_i| = \sum_{a\in T^*} f_i(a)$. Therefore, $(\Lambda_1,\Lambda_2)$ is an equitable subpartition of $\Lambda$, which proves that $G$ is ESP. \qed

\noindent{\bf Remark.} Let us call a graph {\it strong ESP} if every set $\Lambda$ of maximal cliques admits an equitable subpartition $(\Lambda_1,\Lambda_2)$ with $\max\{|\Lambda_1|,|\Lambda_2|\}\le \lceil |\Lambda|/2\rceil$. This proof also proves that (1-)comparability graphs are in fact strong ESP.

\section{Parity graphs}

A graph is called a {\it parity graph} if any two induced paths between the same pair of vertices have the same parity. These are natural extensions of bipartite graphs and they are perfect \cite{sachs}. Cameron and Edmonds \cite{cameron} conjectured that
every parity graph is box-perfect. The objective of this section is to present a proof of this conjecture. 

To establish our result we need a structural characterization of parity graphs. 
Let $H$ be a graph with a stable set $S$ such that all vertices of $S$ have the same set of neighbors. 
Let $B$ be a bipartite graph and let $T$ be a subset of a color class of $B$ with $|T|=|S|$. 
Let $G$ be obtained from the disjoint union of $H$ and $B$ by identifying $S$ with $T$. We call $G$ a {\it bipartite extension} of $H$ by $B$, and we also call the construction of $G$ from $H$ {\it bipartite extension}.

\begin{lemma}[Burlet and Uhry \cite{uhry}]
Every connected parity graph can be constructed from a single vertex by repeatedly duplicating vertices and bipartite extensions.
\end{lemma}

\begin{lemma}\label{lem:twin}
Duplicating a vertex in an ESP graph results in an ESP graph.
\end{lemma}

\vspace{-2mm}
\noindent{\bf Proof.} Let ESP graph $G$ have a vertex $v$. Let $G'$ be obtained by duplicating $v$ and let $v'$ be the new vertex. For any set $\Lambda'$ of maximal cliques of $G'$, we prove that $\Lambda'$ has an equitable subpartition. 

We define $\Lambda$ as follows. 
If $vv'$ is an edge then $\Lambda = \{K-v':K\in \Lambda'\}$; if $vv'$ is not an edge then $\Lambda =\{K:v'\not\in K\in\Lambda\} \uplus \{K-v'+v:v'\in K\in\Lambda'\}$. Note that $\Lambda$ is a multiset of maximal cliques of $G$. Since $G$ is ESP, $\Lambda$ admits an equitable subpartition $(\Lambda_1, \Lambda_2)$. By deleting vertices from  cliques in $\Lambda_1$ and $\Lambda_2$ we may assume that $d_{\Lambda_1} + d_{\Lambda_2} =d_{\Lambda}$ and $\lfloor d_{\Lambda}/2 \rfloor \le d_{\Lambda_i}\le \lceil d_{\Lambda}/2 \rceil$ ($i=1,2$).

If $vv'$ is an edge, let $\Lambda_i' =\{K:v\not\in K\in\Lambda_i\} \uplus \{K+v':v\in K\in\Lambda_i\}$ ($i=1,2$). Then $(\Lambda_1', \Lambda_2')$ is an equitable subpartition of $\Lambda'$ because $d_X(v')=d_X(v)$ holds for $X\in \{\Lambda,\Lambda_1', \Lambda_2'\}$. 

Now suppose $vv'$ is not an edge. Note that $d_{\Lambda} (v) = d_{\Lambda'}(v) + d_{\Lambda'}(v')$. Also we may assume that $d_{\Lambda_1}(v) = \lfloor d_{\Lambda}(v)/2 \rfloor$ and  $d_{\Lambda_2}(v) = \lceil d_{\Lambda}(v)/2 \rceil$. Let 
$$m_1=\lfloor d_{\Lambda'}(v)/2 \rfloor, \ \ m_2=\lceil d_{\Lambda'}(v)/2 \rceil, \ \ m_1'=d_{\Lambda_1}(v)-m_1, \ \ m_2'=d_{\Lambda_2}(v)-m_2.$$
Then 
$$m_1+m_2=d_{\Lambda'}(v), \ \ m_1'+m_2'=d_{\Lambda'}(v'), \ \ \min\{m_1',m_2'\}\ge \lfloor d_{\Lambda'}(v')/2 \rfloor.$$
For $i=1,2$, let $\Lambda_i'$ be obtained from $\Lambda_i$ by turning $m_i'$ cliques $K$ that contain $v$ into $K-v+v'$. Then the above equalities and inequalities imply that $(\Lambda_1', \Lambda_2')$ is an equitable subpartition of $\Lambda'$. \qed

\noindent{\bf Remark.} Clearly, this proof also proves that duplicating a vertex in a strong ESP graph results in a strong ESP graph. 

\begin{theorem}
Parity graphs are ESP.
\end{theorem}

\vspace{-2mm}
\noindent{\bf Proof.} By Lemma \ref{lem:twin}, we only need to show that if $G$ is a bipartite extension of an ESP graph $H$ by a bipartite graph $B=(X,Y,E)$, then $G$ is ESP. Let $X_0\subseteq X$ be the intersection of $H$ and $B$. Let $\Lambda$ be a set of maximal cliques of $G$. Naturally, $\Lambda$ can be partitioned into $\Lambda_H$ and $\Lambda_B$, which are maximal cliques of $H$ and edges of $B$, respectively. Now we find an equitable subpartition $(\Lambda_B',\Lambda_B'')$ of $\Lambda_B$ and an equitable subpartition $(\Lambda_H',\Lambda_H'')$ of $\Lambda_H$ such that $(\Lambda_B'\cup\Lambda_H',\Lambda_B''\cup\Lambda_H'')$ is an equitable subpartition of $\Lambda$. Let $X_0$ be partitioned into $(X_1,X_2)$ such that $X_1$ consists of $x\in X_0$ with both $d_{\Lambda_B}(x)$ and $d_{\Lambda_H}(x)$ odd. Since $(\Lambda_B',\Lambda_B'')$ and $(\Lambda_H',\Lambda_H'')$ are always compatible on vertices in $X_2$, we only need to focus on vertices in $X_1$.

Without loss of generality, let $\Lambda_B=E$. Suppose $B$ has $2t$ vertices of odd degree. Then $E$ can be partitioned into cycles and $t$ paths $P_1,..., P_t$. Let $(\Lambda_B',\Lambda_B'')$ be defined by assigning edges to the two parts alternatively along the cycles and paths. Then $(\Lambda_B',\Lambda_B'')$ is an equitable partition. Note that we have the following freedom in the assignment. Let $x\in X_1$ and let $P_i$ be the path with $x$ as an end. If the other end of $P_i$ is not in $X_1$, then we may choose $d_{\Lambda_B'}(x)$ to be $\lfloor d_{\Lambda_B}(x)/2 \rfloor$ or $\lceil d_{\Lambda_B}(x)/2\rceil$, as we wish (without changing $d_{\Lambda_B'}(z)$ and $d_{\Lambda_B''}(z)$ for any other $z\in X_1$). If the other end of $P_i$ is a vertex $x'$ in $X_1$, then we may assume that $d_{\Lambda_B'}(x)=\lfloor d_{\Lambda_B}(x)/2 \rfloor$ and $d_{\Lambda_B'}(x')=\lceil d_{\Lambda_B}(x')/2\rceil$. Let $(x_1,x_1')$, ..., $(x_k,x_k')$ be these pairs in $X_1$.

Let $H_1$ be obtained from $H$ by deleting $x_1',...,x_k'$ and let $\Lambda_1$ be obtained from $\Lambda_H$ by replacing each $x_i'$ with $x_i$. Note that $d_{\Lambda_1}(x_i) = d_{\Lambda_H}(x_i) + d_{\Lambda_H}(x_i')$ for all $i$, while $d_{\Lambda_1}(v)=d_{\Lambda_H}(v)$ for all other vertices $v$ of $H_1$. Since $H$ is ESP, so is $H_1$. Let $(\Lambda_1', \Lambda_1'')$ be an equitable subpartition of $\Lambda_1$. Without loss of generality, we assume $d_{\Lambda_1'}(x_i)=d_{\Lambda_1''}(x_i) = d_{\Lambda_1}(x_i)/2$ for all $i$. Let $\Lambda_H'$ be obtained from $\Lambda_1'$ by turning $\lfloor d_{\Lambda_H}(x_i')/2 \rfloor$ of its cliques $K$ that contain $x_i$ into $K-x_i+x_i'$ (for every $i$). Then $d_{\Lambda_H'}(x_i)=\lceil d_{\Lambda_H}(x_i)/2\rceil$ and 
$d_{\Lambda_H'}(x_i')=\lfloor d_{\Lambda_H}(x_i')/2 \rfloor$. Let $\Lambda_H''$ be obtained analogously. Now it is straightforward to verify that, the freedom on partition $(\Lambda_B',\Lambda_B'')$ allows us to make adjustments so that  $(\Lambda_B'\cup\Lambda_H',\Lambda_B''\cup\Lambda_H'')$ is an equitable subpartition of $\Lambda$. \qed

\section{Complements of line graphs}

In the rest of this paper we allow some graphs to have loops and parallel edges. We call these {\it multigraphs} and we reserve the word {\it graph} for simple graphs. If a multigraph $H$ is obtained from a graph $H_0$ by adding loops and parallel edges, then $H_0$ is called a {\it simplification} of $H$ and is denoted by $si(H)$.

Let $L(H)$ denote the line graph of a multigraph $H$. Under this circumstance, we always make the following implicit assumptions: \\ 
\indent (i) $H$ has no isolated vertices (deleting an isolated vertex does not affect $L(H)$); \\ 
\indent (ii) $H$ has no loops (replacing a loop with a pendent edge does not affect $L(H)$); \\ 
\indent (iii) $H$ has no distinct vertices $x,y,z$ such that $z$ is the only neighbor of $x$ and the only neighbor \makebox[34pt]{} of $y$ (replacing edges between $y$ and $z$ by edges between $x$ and $z$ does not affect $L(H)$). 

The complement of $L(H)$ will be denoted by $\bar L(H)$. Our results in the next two sections imply a characterization of box-perfect line graphs. The goal of this section is to characterize box-perfect graphs that are complements of line graphs.

\begin{theorem} \label{thm:linebar}
Let $G=\bar L(H)$ be perfect. Then $G$ is box-perfect if and only if $G$ is $\{S_3,\bar S_3^+\}$-free.
\end{theorem}

Our proof of this theorem is divided into a sequence of lemmas. We first determine the structure of $\{S_3,\bar S_3^+\}$-free perfect graphs of the form $\bar L(H)$, and then we confirm that all such graphs are ESP. We will see that some of these graphs are in fact strong ESP.

We need a result of Gallai \cite{gallai} which identifies eight classes and ten individual graphs such that a graph is a comparability graph if and only if it does not contain any of these identified graphs as an induced subgraph. We will use the following immediate consequence of Gallai's theorem. Let $\Gamma$ be the graph obtained from a 6-cycle $v_1v_2v_3v_4v_5v_6v_1$ by adding two edges $v_1v_3$ and $v_1v_5$.

\begin{lemma} \label{lem:inc}
Let $G$ be claw-free and perfect. Then $G$ is an incomparability graph if and only if $G$ does not contain any of $S_3$, $\bar{S}_3$, $\Gamma$, and $C_{2n}\ (n\ge3)$ as an induced subgraph. 
\end{lemma}

Let $K_4^+$ denote the graph obtained from $K_4$ by adding two pendent edges to two of its distinct vertices. Let $K_{2,n}^+$ denote the graph obtained from $K_{2,n}$ ($n\ge3$) by adding a pendent edge to a degree-2 vertex and an edge between the two degree-$n$ vertices.

\begin{lemma}\label{lem:bars3}
Let $\bar L(H)$ be $\{C_5, S_3,\bar S_3^+\}$-free. If $H$ contains $\bar{S}_3$ as a subgraph then $si(H)$ is either $K_4^+$ or a subgraph of $K_{2,n}^+$ for some $n\ge3$.  
\end{lemma}

\vspace{-2mm}
\noindent{\bf Proof.} Since $\bar{S}_3$ is a subgraph of $H$, we assume $V(H)=\{x_1, x_2, x_3, y_1, y_2, y_3, z_1, ..., z_m\}$ such that $x_1x_2x_3$ is a triangle and $x_iy_i\in E(H)$ ($i=1,2,3$). If $m=0$ then it is straightforward to verify the conclusion of the lemma, using the fact that $H$ does not contain $C_5$ as a subgraph. So we assume $m>0$. Let $K_{1,3}^*$ denote the graph obtained from $K_{1,3}$ by subdividing each edge exactly once. Note that $K_{1,3}^*$ is not a subgraph of $H$ since $\bar{L}(K_{1,3}^*)=S_3$. As a result, each $z_i$ is adjacent to none of $y_1,y_2,y_3$, and at most two of $x_1,x_2,x_3$. Furthermore, since $\bar L(H)$ is $\bar S_3^+$-free, the entire neighborhood of each $z_i$ must be a subset of $\{x_1,x_2,x_3\}$ of size one or two (here we also use assumption (i) above). By assumption (iii) above we may assume that each $z_i$ is adjacent to exactly two of $x_1,x_2,x_3$. Since $C_5$ is not a subgraph of $H$, all $z_i$'s must have the same set of neighborhood. Now, since $m>0$, it is straightforward to verify that $si(H)$ is a subgraph of $K_{2,m+3}^+$. \qed

Let $C$ be an even cycle of length $\ge4$. Let $X$ be a stable set of $C$ and let $Y=V(C)-X-N_C(X)$, where $X$ is allowed to be empty. We construct a bipartite graph from $C$ by adding a pendent edge to each vertex in $Y$ and by repeatedly duplicating vertices in $X$. Let $\cal C$ consist of all graphs that can be constructed in this way. 

\begin{lemma} \label{lem:nos3}
Let $L(H)$ be perfect and $\bar S_3$-free. Suppose $H$ is connected and $H$ does not contain $\bar S_3$ as a subgraph.  If $L(H)$ contains an induced $\Gamma$ or $C_{2n}$ $(n\ge3)$, then $si(H)$ is a subgraph of a graph in $\mathcal C\cup\{K_{3,3}\}$.  
\end{lemma}

\vspace{-2mm}
\noindent{\bf Proof.} Suppose $\Gamma$ is an induced subgraph of $L(H)$. Then $H$ has a subgraph with a 4-cycle $x_1x_2x_3x_4$ and two pendent edges $x_1y_1$, $x_2y_2$. Note that $x_1x_3$ and $x_2x_4$ are not edges of $H$ since $\bar S_3$ is not a subgraph of $H$. Let $z_1,...,z_m$ be the remaining vertices of $H$. If $m=0$, then either $si(H)$ is a subgraph of $K_{3,3}$ or $H$ contains a 5-cycle. So we assume $m>0$. Like in the proof of the last lemma, since $C_5$ and $K_{1,3}^*$ are not subgraphs of $H$, for each $i$ we must have $N_H(z_i)=\{x_1,x_3\}$ or $\{x_2,x_4\}$, or $\{x_j\}$ for some $j$. In addition, $N_H(y_i)\subseteq\{x_i,x_{i+2}\}$ ($i=1,2$) and $|N_H(y_1)\cup N_H(y_2)|\le 3$. Now, since $H$ does not contain $K_{1,3}^*$, it is routine to check that $si(H)$ is a subgraph of a graph in $\cal C$.

Next, suppose $L(H)$ is $\Gamma$-free. Then $H$ contains a $2n$-cycle $x_1x_2... x_{2n}$ ($n\ge3$). Note that this cycle has no chord (otherwise $L(H)$ contains an induced $\Gamma$, $\bar S_3$, or $C_{2k+1}$ with $k\ge2$). Let $z_1,...,z_m$ be the remaining vertices of $H$. Using the same argument we used in the last paragraph it is straightforward to show that each $N_H(z_i)$ is $\{x_j\}$ or $\{x_j,x_{j+2}\}$ for some $j$ (where $x_{2n+t}$ is $x_t$). In addition, if $N_H(z_i)=\{x_j,x_{j+2}\}$ then $N_H(x_{j+1})=N_H(z_i)$. Therefore, $si(H)$ is a subgraph of a graph in $\cal C$. \qed

\begin{lemma} \label{lem:ab}
Suppose $G$ has a vertex $u$ such that $G-u$ is bipartite and $G-N(u)$ is edge-less. Then $G$ is totally unimodular.  
\end{lemma}

\vspace{-2mm}
\noindent{\bf Proof.} By Theorem 19.3 of \cite{schrijver}, we only need to show that each set $\Lambda$ of maximal cliques admits an {\it equitable} partition $(\Lambda_1,\Lambda_2)$, meaning that $\min\{d_{\Lambda_1}(v), d_{\Lambda_2}(v)\} \ge \lfloor d_{\Lambda}(v)\rfloor$, for all $v\in V(G)$. Suppose to the contrary that some $\Lambda$ does not admit such a partition. We choose $\Lambda$ with $|\Lambda|$ as small as possible.  

Let $A,B,C,D$ be a partition of $V(G)-u$ such that $A\cup C$, $B\cup D$ are stable and $N(u)=B\cup C$. Let $G'$ be the subgraph of $G$ formed by edges in $K-u$, over all $K\in \Lambda$. We claim that $G'$ is a forest. Suppose $G'$ has a cycle $x_1x_2...x_n$. Note that for each $i$, exactly one of $x_ix_{i+1}$ and $ux_ix_{i+1}$ is a clique in $\Lambda$. Let $\Lambda'$ be the rest cliques in $\Lambda$. By the minimality of $|\Lambda|$, $\Lambda'$ admits an equitable partition $(\Lambda_1',\Lambda_2')$. Let us extend $\Lambda_j'$ ($j=1,2$) to $\Lambda_j$ by including $x_ix_{i+1}$ or $ux_ix_{i+1}$ (whichever belongs to $\Lambda$) for all $i$ with $i-j$ even. Then it is easy to see that $(\Lambda_1,\Lambda_2)$ is an equitable partition of $\Lambda$. This contradicts the choice of $\Lambda$ and thus the claim is proved. The same argument also shows that $G'$ has no maximal path with two ends both in $A\cup B$ or both in $C\cup D$. Thus all components of $G'$ are paths with one end in $A\cup B$ and one end in $C\cup D$. If $G'$ has only one path then the same argument still works. If $G'$ has two or more paths then we can take any two of them and treat their union as a cycle and again apply the same argument. \qed

Recall that a graph $G$ is {\it strong ESP} if every set $\Lambda$ of maximal cliques of $G$ admits an equitable subpartition $(\Lambda_1,\Lambda_2)$ with $\max\{|\Lambda_1|,|\Lambda_2|\}\le \lceil |\Lambda|/2\rceil$. The next lemma follows immediately from this definition.

\begin{lemma} \label{lem:strongesp}
(1) If $G$ is strong ESP then so are all its induced subgraphs. \\ 
\indent (2) Let $G_1,G_2$ be strong ESP and let $G$ be obtained from the disjoint union of $G_1,G_2$ by adding all edges between them. Then $G$ is also strong ESP. 
\end{lemma}

In a (loopless) multigraph $G$, the {\it degree} of a vertex $v$, denoted $d_G(v)$, is the number of edges incident with $v$. The next is the key step for proving Theorem \ref{thm:linebar}.

\begin{lemma} \label{lem:k33}
For every $H\in\mathcal C\cup\{K_{3,3}\}$, $\bar L(H)$ is strong ESP.
\end{lemma}

\vspace{-2mm}
\noindent{\bf Proof.} For each $\mu\in \mathbb Z_+^{E(H)}$, let $\mu H$ denote the multigraph with vertex set $V(H)$ such that the number of edges  between any two vertices $x,y$ is zero (if $xy\not\in E(H)$) or $\mu(xy)$ (if $xy\in E(H)$). Note that $\mu H$ is bipartite since $H$ is bipartite. Let $\Delta(\mu)$ denote the maximum degree of $\mu H$. By Konig's edge-coloring theorem, $E(\mu H)$ is the union of $k$ matchings if and only if $k\ge \Delta(\mu)$. Because of this theorem and the one-to-one correspondence between cliques of $\bar L(H)$ and matchings of $H$, to prove the lemma it is enough for us to show that 

($*$) \ for any $\mu\in \mathbb Z_+^{E(H)}$ there exist $\mu_1,\mu_2\in \mathbb Z_+^{E(H)}$ such that $\mu_1+\mu_2=\mu$, $\mu_i\ge\lfloor \mu/2 \rfloor$ ($i=1,2$), \\ 
\makebox[35pt]{} $\Delta(\mu_1) \le \lceil \Delta(\mu)/2\rceil$, and $\Delta(\mu_2) \le \lfloor \Delta(\mu)/2\rfloor$.

In the following we construct a partition $(E_1,E_2)$ of $E(\mu H)$ such that the multiplicity functions $\mu_i$ of $E_i$ ($i=1,2$) satisfies ($*$). This partition will be constructed in several steps. In the process we determine a partition $(E_1,E_2,E_3)$ of $E(\mu H)$, where we begin with $(E_1,E_2,E_3)=(\emptyset,\emptyset, E(\mu(H))$ and we keep moving edges from $E_3$ to $E_1,E_2$ until $E_3$ becomes empty. For $i=1,2,3$, let $H_i$ denote the subgraph of $\mu H$ formed by edges in $E_i$. 

First, for each edge $e=xy$ of $H$, among all $\mu(e)$ edges of $E_3$ that are between $x$ and $y$, we move $\lfloor \mu(e)/2 \rfloor$ of them to $E_1$ and $\lfloor \mu(e)/2 \rfloor$ of them to $E_2$. At the end of this process, $H_3$ becomes a simple graph. It follows that $\mu_i\ge\lfloor \mu/2 \rfloor$ ($i=1,2$) and this inequality will be satisfied no matter how edges of $H_3$ are moved to $E_1$ and $E_2$ in later steps.

If $H_3$ has a cycle $C$, since $H$ is bipartite, $E(C)$ can be partitioned into two matchings $M_1,M_2$. We move $M_i$ from $E_3$ to $E_i$ ($i=1,2$). We repeat this process until $H_3$ become a forest. At this point, $H_1$ and $H_2$ have the same degree on every vertex.

Let $S=\{v: d_{\mu H}(v)=\Delta(\mu)\}$. Suppose $H_3$ has a leaf $v$ that is not in $S$. Let $P$ be a maximal path of $H_3$ starting from $v$. Let $E(P)$ be partitioned into two matchings $M_1,M_2$, where we assume the edge of $P$ that is incident with the other end $u$ of $P$ belongs to $M_1$. Then we move $M_i$ from $E_3$ to $E_i$ ($i=1,2$). After this change, $d_{H_1}(u)=\lceil d_{\mu H}(u)/2 \rceil \le \lceil \Delta(\mu)/2\rceil$,  $d_{H_2}(u)=\lfloor d_{\mu H}(u)/2 \rfloor \le \lfloor \Delta(\mu)/2\rfloor$, and $d_{H_i}(v)\le \lceil d_{\mu H}(v)/2\rceil \le \lfloor \Delta(\mu)/2\rfloor$ ($i=1,2$). In addition, $d_{H_1}(w)= d_{H_2}(w)$ for all $w\ne u,v$, and $d_{H_i}(u),d_{H_i}(v)$ will remain unchanged in the remaining process. By repeating this process we may assume that all leaves of $H_3$ are in $S$. As a consequence, $\Delta(\mu)$ is odd. Note that the same argument works if $H_3$ has a maximal path with an odd number of edges. Thus we further assume that in every component of $H_3$, all leaves are in the same color class (of any 2-coloring of $H_3$).

We first consider the case $H=K_{3,3}$. We claim that each component of $H_3$ is a path. Suppose a component $H_3'$ of $H_3$ is not a path. Then $H_3'$ has at least three leaves. Since all these leaves are in the same color class, $H_3'$ must have exactly three leaves $z_1,z_2,z_3$ and they form a color class of $H$. Consequently, $H_3'=H_3=K_{1,3}$. Moreover, in the previous steps of reducing $H_3$, no path was ever deleted because otherwise $H_3$ would be a subgraph of $K_{2,2}$. It follows that $d_{\mu H}(v^*)$ is even, where $v^*\in V(H)-V(H_3)$. However, the fact $z_1,z_2,z_3\in S$ implies that $\mu H$ is $\Delta(\mu)$-regular, and thus $d_{\mu H}(v^*)=\Delta(\mu)$ is odd. This contradiction proves our claim. Now, since each non-leaf $v$ of $H_3$ has degree two, its degree in $\mu H$ is even and thus $v\not\in S$. It follows that moving all edges of $E_3$ to $E_1$ results in the required partition. 

Next suppose $H\in\cal C$. Let $H_3'$ be a component of $H_3$. Then $H_3'$ is a {\it caterpillar} since $K_{1,3}^*$ is not a subgraph of $H$. Therefore, $H_3'$ has a path $x_1x_2...x_{2k+1}$ such that every leaf of $H_3'$ is adjacent to some $x_{2i+1}$. We assume that $H_3'$ is not a path because otherwise we may move the entire path from $E_3$ to $E_1$. We make two observations before we continue. First, $d_H(v)>1$ holds for every leaf $v$ of $H_3'$, because otherwise the only edge of $H$ that is incident with $v$ would be the only edge of $H_3'$ (as $v\in S$). Second, if $u,v\in V(H_3')$ are of degree-2 in $H$ and are contained in a 4-cycle $uxvy$ of $H$, then at most one of $u,v$ is in $S$. This is because otherwise $\mu(ux)=\mu(vy)$, $\mu(uy)=\mu(vx)$, and both $x,y\in S$, which implies that $H_3'$ is a subgraph of the 4-cycle $uxvy$. It follows from these two observations and the construction of graphs in $\cal C$ that each $x_{2i+1}$ is adjacent to at most two leaves of $H_3'$. For the same reasons, there must exist $i_0\in\{0,1,...,k\}$ such that $d_{H_3'}(x_{2i_0+1})=2$.

For $i=1,2,3,4$, let $V_i=\{v: d_{H_3'}(v)=i\}$. Note that $V_2\cup V_3\cup V_4 =\{x_1,...,x_{2k+1}\}$. Let $M$ be the matching $\{x_{2i-1}x_{2i}:i=1,...,i_0\} \cup \{x_{2i}x_{2i+1}:i_0+1,...,k\}$. From $H_3'$ we move $M$ to $E_2$ and the rest of $E(H_3')$ to $E_1$. Now we verify that, after this change, $d_{H_1}(v)\le \lceil \Delta(\mu)/2 \rceil$ and $d_{H_2}(v) \le \lfloor \Delta(\mu)/2\rfloor$ hold for all $v\in V_1\cup V_2\cup V_3\cup V_4$. For each $v\in V_1$ it is easy to see that in fact $d_{H_1}(v)=\lceil \Delta(\mu)/2 \rceil$ and $d_{H_2}(v) = \lfloor \Delta(\mu)/2\rfloor$. For each even $i$, we have $x_i\in V_2$ and $d_{H_1}(x_i) = d_{H_2}(x_i) = d_{\mu H}(x_i)/2\le \lfloor \Delta(\mu)/2\rfloor$. For each odd $i$ we consider two cases. If $x_i\in V_3$ then $d_{H_1}(x_i)= (d_{\mu H}(x_i)+1)/2\le \lceil \Delta(\mu)/2 \rceil$ and $d_{H_2}(x_i)= (d_{\mu H}(x_i)-1)/2\le \lfloor \Delta(\mu)/2\rfloor$. If $x_i\in V_2\cup V_4$ then $d_{H_1}(x_i)\le (d_{\mu H}(x_i)+2)/2\le \lceil \Delta(\mu)/2 \rceil$ and $d_{H_2}(x_i)\le d_{\mu H}(x_i)/2\le \lfloor \Delta(\mu)/2\rfloor$. Therefore, we may apply this split to all components of $H_3$ and create the required partition $E_1,E_2$. \qed

\bigskip
\noindent{\bf Proof of Theorem \ref{thm:linebar}.} The forward implication is obvious so we only show that $G=\bar L(H)$ is ESP when $G$ is perfect and $\{S_3, \bar{S}_3^+\}$-free. 

Suppose $L(H)$ contains an induced $S_3$. Then $H$ contains $\bar S_3$ as a subgraph. By Lemma \ref{lem:bars3}, $si(H)$ is either $K_4^+$ or a subgraph of $K_{2,n}^+$ for some $n\ge3$. In both cases, it is straightforward to verify that $\bar L(si(H))$ satisfies the assumptions in Lemma \ref{lem:ab}. So $\bar L(si(H))$ is totally unimodular and thus is also ESP. By Lemma \ref{lem:twin}, $\bar L(H)$ is ESP.

Now suppose $L(H)$ is $S_3$-free. We claim that $\bar L(H')$ is strong ESP for every component $H'$ of $H$. If $\bar L(H')$ is a comparability graph, then the claim follows immediately from the Remark at the end of Section 4. So we assume that $\bar L(H')$ is not a comparability graph. By Lemma \ref{lem:inc}, $L(H)$ contains an induced $\Gamma$ or $C_{2n}$ ($n\ge3$). This implies, by Lemma \ref{lem:nos3}, that $si(H')$ is a subgraph of a graph in $\mathcal C\cup\{K_{3,3}\}$. Then the claim follows from Lemma \ref{lem:k33}, Lemma \ref{lem:strongesp}(1), and the Remark of Lemma \ref{lem:twin}. Finally, this claim and Lemma \ref{lem:strongesp}(2) imply that $\bar L(H)$ is ESP. \qed

\section{Trigraphs}

Our next objective is to characterize claw-free box-perfect graphs. To accomplish this goal, we will need a result of Chudnovsky and Plumettaz \cite{CP} on the structure of claw-free perfect graphs. The purpose of this section is to explain their result, which requires many definitions. 

A \emph{trigraph} $G$ consists of a finite set $V$ of {\it vertices} and an {\it adjacency function} $\theta: \binom{V}{2} \rightarrow \{1, 0, -1\}$ such that $\{uv: \theta(uv)=0\}$ is a matching. 
Two distinct vertices $u$ and $v$ of $G$ are \emph{strongly adjacent} if $\theta(uv)=1$, \emph{strongly antiadjacent} if $\theta(uv)=-1$, and \emph{semiadjacent} if $\theta(uv)=0$. We call $u$, $v$ \emph{adjacent} if $\theta(uv)\geq 0$, and \emph{antiadjacent} if $\theta(uv)\leq 0$. Note that every graph can be considered as a trigraph with $\{uv: \theta(uv)=0\}=\emptyset$. In other words, graphs are exactly trigraphs with no semiadjacent pairs. The result of Chudnovsky and Plumettaz is in fact about trigraphs.

For any trigraph $G =(V,\theta)$, let $G^{\ge0}$ denote the graph $(V,\{uv:\theta(uv)\ge0\})$. Conversely, for any graph $G=(V,E)$, let $tri(G)$ denote the set of all trigraphs $(V,\theta)$ such that for any distinct $u,v\in V$, $\theta(uv)\ge0$ if $uv\in E$ and $\theta(uv)\le 0$ if $uv\not\in E$. 

Let $G=(V,\theta)$ be a trigraph. 
We call $G$ {\it connected} if $G^{\ge0}$ is connected. 
For each $v\in V$, let $N_G(v)=N_{G^{\ge0}}(v)$. We often write $N(v)$ for $N_G(v)$ if the dependency on $G$ is clear. 
For any $X\subseteq V$, let $G|X$ be the trigraph such that its vertex set is $X$ and its adjacency function is the restriction of $\theta$ to $\binom{X}{2}$. If a trigraph $H$ is  isomorphic to $G|X$ for some $X\subseteq V$, then we call $H$  a {\it subtrigraph} of $G$ and we say that $G$ {\it contains} $H$.

A trigraph is a {\it hole} if it belongs to $tri(C_n)$ for some $n\ge4$. A trigraph $(V,\theta)$ is an {\it antihole} if $(V,-\theta)$ is a hole. A hole or antihole is {\it odd} if its number of vertices is odd. A trigraph is {\it Berge} if it contains neither odd hole nor odd antihole. A trigraph is a {\it claw} if it belongs to $tri(K_{1,3})$. 
A trigraph is {\it claw-free} if it does not contain any claw. In general, if $\cal H$ is a set of trigraphs, then a trigraph is {\it $\cal H$-free} if it does not contain any trigraph in $\cal H$. The result of Chudnovsky and Plumettaz characterizes \{claw, holes, antiholes\}-free trigraphs, that is, claw-free Berge trigraphs. To describe the resulting structure we need more definitions.

Let $G=(V,\theta)$ be a trigraph. For any two disjoint $X, Y\subseteq V$, we say that $X$ is \emph{complete} (resp. \emph{strongly complete}, \emph{anticomplete}, \emph{strongly anticomplete}) to $Y$ if every $x\in X$ and every $y\in Y$ are adjacent (resp. strongly adjacent, antiadjacent, strongly antiadjacent). A \emph{clique} (resp. \emph{strong clique}) of $G$ is a set $C\subseteq V$ such that any two distinct vertices of $C$ are adjacent (resp. strongly adjacent). A \emph{stable set} (resp. \emph{strong stable set}) of $G$ is a set $S\subseteq V$ such that any two distinct vertices of $S$ are antiadjacent (resp. strongly antiadjacent).

A trigraph $H$ is a \emph{thickening} of a trigraph $G$ if $V(H)$ admits a partition $(X_v:v\in V(G))$ such that \\ 
\indent$\bullet$ \ if $v\in V(G)$ then $X_v\ne\emptyset$ is a strong clique of $H$;\\ 
\indent$\bullet$ \ if $u, v \in V(G)$ are strongly adjacent in $G$ then $X_u$ is strongly complete to $X_v$ in $H$; \\ 
\indent$\bullet$ \ if $u, v \in V(G)$ are strongly antiadjacent, then $X_u$ is strongly anticomplete to $X_v$ in $H$;\\ 
\indent$\bullet$ \ if $u, v \in V(G)$ are semiadjacent, then $X_u$ is neither strongly complete nor strongly \\ \makebox[9mm]{} anticomplete to $X_v$ in $H$.

\begin{figure}[ht]
\centerline{\includegraphics[scale=0.5]{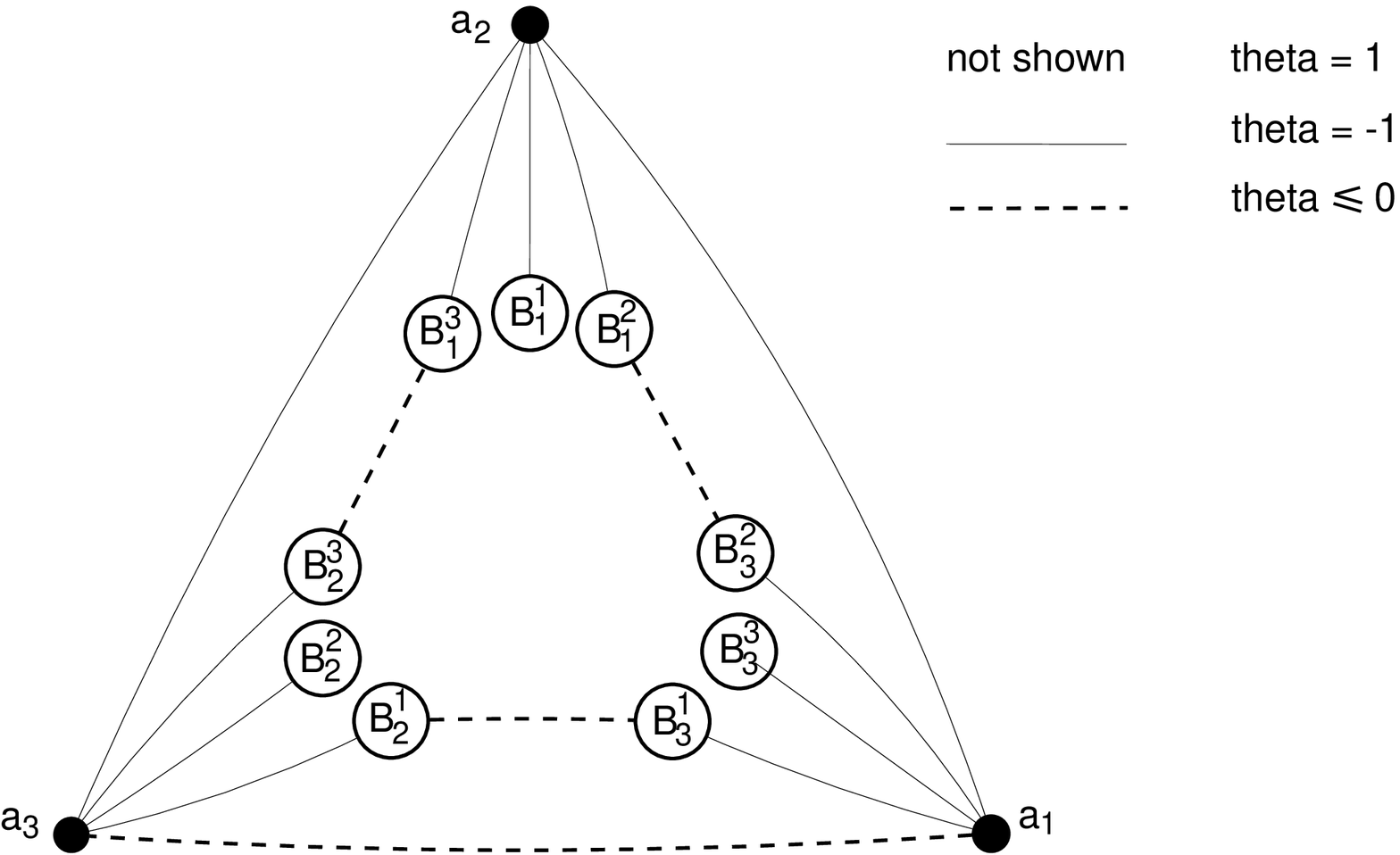}}
%\centerline{\includegraphics[scale=0.5]{Figures/C.eps}}
\caption{Trigraphs in $\cal C$}
\label{fig:c}
\end{figure}

Let ${\cal C}$ be the class of all trigraphs illustrated in Figure \ref{fig:c}, where \\ 
\indent $\bullet$ \ $|B^j_i| \leq 1$  for all $i,j\in\{1,2,3\}$\\ 
\indent $\bullet$ \ $|B_2^1\cup B_3^1|$, $|B_1^2\cup B_3^2|$, $|B_1^3\cup B_2^3|\in \{0,2\}$ \\ 
\indent $\bullet$ \ if $\theta(a_1a_3)=0$ then $B_2^1\cup B_3^1 =\emptyset$\\ 
\indent $\bullet$ \ there exists $x_i\in B_i^1\cup B_i^2\cup B_i^3$ for $i=1,2,3$, such that $\{x_1,x_2,x_3\}$ is a clique.

It turns out that there are two kinds of claw-free Berge trigraphs. The first are thickenings of trigraphs in $\cal C$. The second are constructed (in a way like constructing line graphs) from certain basic trigraphs. In the following, we first define the building blocks and then describe the construction.

Let $G$ have three vertices $v, z_1, z_2$ such that $\theta(vz_1)=\theta(vz_2)=1$ and $\theta(z_1z_2)=-1$. Then the pair $(G, \{z_1, z_2\})$ is a {\it spot}. Let $G$ have four vertices $v_1,v_2,z_1,z_2$ such that $\theta(v_1z_1)=\theta(v_2z_2)=1$, $\theta(v_1v_2)=0$, $\theta(z_1z_2)=\theta(z_1v_2)=\theta(z_2v_1)=-1$. Then the pair $(G,\{z_1, z_2\})$ is a {\it spring}. 

A trigraph is a {\it linear interval} if its vertices can be ordered as $v_1,...,v_n$ such that if $i < j < k$ and $\theta(v_iv_k)\ge 0$ then $\theta(v_iv_j) = \theta(v_jv_k)=1$. Let $G$ be such a trigraph with $n\ge4$. We call $(G, \{v_1, v_n\})$ a {\it linear interval stripe} if: $v_1$ and $v_n$ are strongly antiadjacent, $v_i$ and $v_{i+1}$ are adjacent for every $i\in\{1,...,n-1\}$, no vertex is complete to $\{v_1, v_n\}$, and no vertex is semiadjacent to $v_1$ or $v_n$.

Let $(G,\{p,q\})$ be a spring or a linear interval strip. Let $H$ be a thickening of $G$ and let $X_v$ ($v\in V(G)$) be the corresponding sets. If $|X_p|=|X_q|=1$, then $(H,X_p\cup X_q)$ is called a {\it thickening} of  $(G,\{p,q\})$.

Let $\mathcal C'$ be the class of all pairs $(H, \{z\})$ such that $H$ is a thickening of a trigraph $G\in\cal C$ and $z\in X_{a_i}$ for some $i\in \{1,2,3\}$  for which $B^{i+2}_{i+1}\cup B^{i+2}_{i}=\emptyset$ and $N(z)\cap (X_{a_{i+1}}\cup X_{a_{i+2}})=\emptyset$ (here we use the notation from the definitions of $\cal C$ and thickening).

A {\it signed graph} $(G,s)$ consists of a multigraph $G=(V,E)$ and a function $s: E\rightarrow \{0, 1\}$. If $\sum_{e\in E(C)}s(e)$ is even for all cycles $C$ of $G$, then $(G, s)$ is an \emph{evenly signed graph}. In the following we define another three classes of signed graphs. 
For any $F\subseteq E$, let $G[F]= (V, F)$.

Let $\mathcal F_1$ be the class of loopless signed graphs $(G,s)$ such that $si(G)=K_4$ and $s \equiv 1$.  
Let $\mathcal F_2$ be the class of loopless signed graphs $(G,s)$ such that $si(G)$ is obtained from $K_{2,n}$ ($n\ge1$) by adding an edge $e^*$ between its two degree-$n$ vertices, and edges in $\{e:s(e) = 0\}$ are all parallel to $e^*$ (while $s(e^*)=1$). 
We remark that our $\mathcal F_1$ is $\mathcal F_2$ of \cite{CP} and our  $\mathcal F_2$ is $\mathcal F_1\cup \mathcal F_3$ of \cite{CP}

In a connected multigraph $G$ with $E(G)\ne\emptyset$, a subgraph $B$ is a {\it block} of $G$ if $B$ is a loop or $B$ is maximal with the property that $B$ is loopless and $si(B)$ is a block of $si(G)$. A signed graph $(G, s)$ is called an \emph{even structure} if $E(G)\ne\emptyset$ and for all blocks $B$ of $G$, $(B, s|_{E(B)})$ is a member of ${\cal F}_1\cup {\cal F}_2$ or an evenly signed graph or a loop.

Now we describe how the pieces defined above can be put together. A trigraph $G=(V,\theta)$ is called an \emph{evenly structured linear interval join} if it can be constructed in the following manner: \\ 
\indent $\bullet$ Let $(H,s)$ be an even structure. \\ 
\indent $\bullet$ For each edge $e\in E(H)$, let $Z_e\subseteq V(H)$ be the set of ends of $e$ (so $|Z_e|=1$ or 2). \\ \makebox[23pt]{} Let $S_e=(G_e,Z_e)$ such that $G_e$ is a trigraph with $V(G_e)\cap V(H)=Z_e$ and \\ 
\indent\indent $*$ if $e$ is not on any cycle then $S_e$ is a spot or a thickening of a linear interval stripe, \\ 
\indent\indent $*$ if $e$ is on a cycle of length $>1$ and $s(e)=0$ then $S_e$ is a thickening of a spring,\\ 
\indent\indent $*$ if $e$ is on a cycle of length $>1$ and $s(e)=1$ then $S_e$ is a spot, \\ 
\indent\indent $*$ if $e$ is a loop then $S_e\in {\cal C}'$. \\ 
\indent $\bullet$ For all distinct $e,f\in E(H)$, $V(G_e)\cap V(G_f)\subseteq Z_e\cap Z_f$. \\ 
\indent $\bullet$ Let $V=\cup_{e\in E(H)} V(G_e)\backslash Z_e$ and let $\theta$ be given by: for any $u,v\in V$ \\
\indent\indent $*$ if $u,v\in V(G_e)\backslash Z_e$ for some $e\in E(H)$ then $\theta(uv)=\theta_{G_e}(uv)$ \\ 
\indent\indent $*$ if $u\in\! N_{G_e}(x)$ and $v\in\! N_{G_f}(x)$ for distinct $e,f\in\! E(H)$ with a common end $x$, then $\theta(uv)=\!1$ \\ 
\indent\indent $*$ in all other cases, $\theta(uv)=-1$. \\ 
\indent $\bullet$ We will write $G=\Omega(H,s,\{S_e:e\in E(H)\})$.

\begin{theorem}[Chudnovsky and Plumettaz \cite{CP}] \label{thm:cp}
A connected trigraph is claw-free and Berge if and only if it is a thickening of a trigraph in $\cal C$ or an evenly structured linear interval join.
\end{theorem}

In the following we produce a different formulation of this result. A vertex $x$ of a trigraph is {\it simplicial} if $N(x)\ne\emptyset$ and $\{x\}\cup N(x)$ is a strong clique. For $i=1,2$, let $G_i=(V_i,\theta_i)$ be a trigraph with a simplicial vertex $x_i$ and with $|V_i|\ge3$. The {\it simplicial sum} of $G_1,G_2$ (over $x_1,x_2$) is the trigraph $G=(V,\theta)$ such that $V=(V_1-x_1)\cup (V_2-x_2)$ and, for all distinct $v_1$, $v_2\in V$, \\ 
\indent $\bullet$ $\theta(v_1v_2)=\theta_i(v_1v_2)$ if $\{v_1,v_2\}\subseteq V_i$ for some $i=1,2$ \\ 
\indent $\bullet$ $\theta(v_1v_2)=1$ if $v_i\in N_{G_i}(x_i)$ for both $i=1,2$ \\ 
\indent $\bullet$ $\theta(v_1v_2)=-1$ if otherwise. \\ 
We point out that both $G_1$ and $G_2$ are contained in $G$. Moreover, using the language of \cite{CP}, $G$ admits either a 1-join or a homogeneous set of size $\ge2$. 

\begin{lemma}\label{lem:sum}
Let $G$ be a simplicial sum of $G_1,G_2$. Then $G$ is claw-free if and only if both $G_1,G_2$ are; and $G$ is Berge if and only if both $G_1,G_2$ are. 
\end{lemma}

We omit the proof since it is straightforward. This lemma suggests that we can characterize claw-free Berge trigraphs by determining all such trigraphs that are not simplicial sums. In the following we describe these trigraphs.

\begin{figure}[ht]
%\centerline{\includegraphics[scale=0.5]{Figures/J.eps}}
\centerline{\includegraphics[scale=0.5]{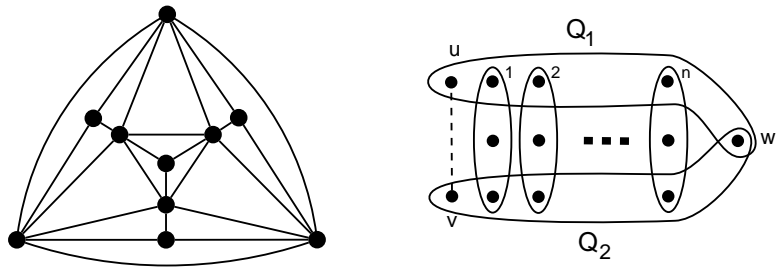}}
\caption{$J_1$ and $J_2$}
\label{fig:j}
\end{figure}

Let $\cal I$ be the class of linear interval trigraphs. 
Let $\cal L$ be the class of trigraphs 
$G$ such that $G^{\ge0}$ is the line graph of a bipartite multigraph and every {\it triangle} (a clique of size 3) of $G$ is a strong clique. Let $J_1$ be the first graph in Figure \ref{fig:j}. We consider $J_1$ as a trigraph with no semiadjacent pairs. Let $\mathcal J_1$ consists of trigraphs obtained from $J_1$ by deleting $k$ of its cubic vertices $(0\le k\le 4)$. Let $J_2(n)$ be the second trigraph in Figure \ref{fig:j}, where $Q_1,Q_2$, and all vertical triples are strong cliques, $\theta(uv)$ could be 0, 1, or $-1$, and all other pairs are strongly antiadjacent. Note that $J_2(0)\in\cal I$. Let $\mathcal J_2$ consist of trigraphs of the form $J_2(n)-X$ for all $n\ge 1$ and all $X\subseteq \{u,v\}$. Let $\mathcal J=\mathcal J_1\cup \mathcal J_2$.

\begin{theorem}\label{thm:sum}
A connected trigraph is claw-free and Berge if and only if it is obtained by simplicial summing thickenings of trigraphs in $\cal C\cup L\cup I\cup J$.
\end{theorem}

We need a few lemmas in order to prove this theorem. A {\it $1$-separation} of a multigraph $H$ is a pair $(H_1,H_2)$ of edge-disjoint proper subgraphs of $H$ such that $H_1\cup H_2=H$ and $|V(H_1)\cap V(H_2)|=1$. Suppose $G=\Omega(H,s,\{S_e\})$. Then a 1-separation $(H_1,H_2)$ of $H$ is called {\it trivial} if there exists $i\in\{1,2\}$ such that $H_i=K_2$ and $S_f$ is a spot, where $f$ is the only edge of $H_i$.

\begin{lemma} \label{lem:trivial}
Suppose $G=\Omega(H,s,\{S_e\})$ and suppose $H$ has a nontrivial 1-separation $(H_1,H_2)$. Then $G$ is a simplicial sum of two trigraphs. 
\end{lemma}

\vspace{-2mm}
\noindent\textbf{Proof.} Let $x$ be the common vertex of $H_1,H_2$. For $i=1, 2$, let $H_i'$ be obtained from $H_i$ by adding a new vertex $x_i$ and a new edge $xx_i$. Let $s_i$ be the signing of $H_i'$ which agrees with $s$ on $H_i$, and $s_i(xx_i)=1$. Since all blocks of $H_i'$ (other than $xx_i$) are blocks of $H$, $(H_i', s_i)$ is an even structure. Let $S_{xx_i}$ be a spot and let $G_i=\Omega(H_i',s_i,\{S_e:e\in E(H_i')\})$. Since separation $(H_1,H_2)$ is nontrivial, $G_i$ must have $\ge 3$ vertices. Now it is straightforward to verify that $x_i$ is a simplicial vertex of $G_i$ ($i=1,2$) and $G$ is the simplicial sum of $G_1$ and $G_2$ over $x_1$ and $x_2$. \qed

\begin{lemma}\label{lem:thick}
Let $H$ be a thickening of $G$. \\ 
%\indent (i) If $F$ is a thickening of $H$ then $F$ is a thickening of $G$. \\ 
\indent (i) $H$ is claw-free if and only if $G$ is claw-free. \\ 
\indent (ii) $H$ is Berge if and only if $G$ is Berge.   
\end{lemma}

\vspace{-2mm}
\noindent{\bf Proof.} Part (ii) is (6.4) of \cite{CP} and part (i) is easy to verify, as pointed out in \cite{claw5}.\qed

A trigraph $G$ is {\it quasi-line} if $N(v)$ is the union of two strong cliques for every $v\in V(G)$. It is easy to see that if $G$ is quasi-line then $G$ is claw-free. A trigraph $G$ is {\it cobipartite} if $V(G)$ is the union of two strong cliques. Clearly, if $G$ cobipartite then $G$ is quasi-line and thus is claw-free. It is also clear that every connected cobipartite trigraph with $\ge 2$ vertices is a thickening of a two-vertex trigraph. Thus every cobipartite trigraph is Berge.

\bigskip
\noindent {\bf Proof of Theorem \ref{thm:sum}.} To prove the backward implication, by Lemma \ref{lem:sum} and Lemma \ref{lem:thick}, we only need to consider trigraphs $G\in \cal C\cup L\cup I\cup J$. If $G\in\cal C$ then the result follows from Theorem \ref{thm:cp}. If $G\in\cal I$ then $G$ is claw-free \cite{claw5} and Berge \cite{CP}. If $G\in\cal L\cup J$ then $G$ is quasi-line and thus $G$ is claw-free. If $G\in \cal J$, then deleting simplicial vertices from $G$ results in a cobipartite trigraph, which implies that $G$ is Berge. Finally, assume $G\in \cal L$ and $G^{\ge0}=L(B)$ is the line graph of a bipartite multigraph $B$. We need to show that $G$ is Berge. Since no semiadjacent pairs are contained in a triangle, every hole of $G$ must come from a cycle of $B$ and thus $G$ contains no odd holes. If $G$ has an antihole $v_1v_2...v_nv_1$ with $n\ge 7$, then we consider the restriction of $G$ on $v_1,...,v_6$. If $\theta(v_iv_{i+1})=-1$ for all $i=1,...,5$, then the graph $X$ formed by $\{v_iv_j:\theta(v_iv_j)\ge0\}$ would be the complement of a path on six vertices, which is one of the minimal non-line-graphs. This is impossible since $X$ is an induced subgraph of $L(B)$. So $\theta(v_iv_{i+1})=0$ holds for some $i$, which makes $v_i,v_{i+1},x_k$ a triangle for some $k$. This contradiction (two semiadjacent vertices are contained in a triangle) shows that $G$ contains no antihole of length $\ge7$. Thus $G$ is Berge, which completes the proof of the backward direction.

To prove the forward implication, by Theorem \ref{thm:cp}, we assume $G=\Omega(H,s,\{S_e\})$. Since $G$ is connected, $H$ is connected as well. By Lemma \ref{lem:trivial}, we also assume that all 1-separations of $H$ are trivial. Let $U$ be the set of all degree-one vertices $u$ of $H$ for which if $e$ is the only edge incident with $u$ then $S_e$ is a spot. We assume $V(H)\ne U$ because otherwise $H=K_2$ and $G=K_1$ and thus the result holds.  
Let $H_0=H-U$. Note that $H_0$ is connected, as $H$ is connected. Moreover, by its construction, $H_0$ dose not have a 1-separation. Thus either $H_0=K_1$ or $H_0$ is a block of $H$.

Suppose $H_0$ is $K_1$ or $K_2$. It follows that $H$ is a tree with 1, 2, or 3 edges. Moreover, $S_e$ is a thickening of a linear interval strip for at most one $e$, and every other $S_e$ is a spot. In all cases, it is routine to check that $G$ is a thickening of a trigraph in $\cal I$.

Suppose $H_0$ is a loop $e$. Let $S_e=(G_e,\{z\})\in\cal C'$ and let $G_e$ be a thickening of $C\in\cal C$. If $H$ has $\ge2$ edges then $H$ consists of $e$ and a pendent edge $f$ with $S_f$ a spot. It follows that $G=G_e$, which is a thickening of a trigraph in $\cal C$. So $e$ is the only edge of $H$ and $G=G_e-z$. If $z$ is not the only vertex of $X_{a_i}$ (here we use the notation in the definition of $\cal C'$) then $G$ is also a thickening of $C$. If $z$ is the unique vertex of $X_{a_i}$ then $G$ is cobipartite. In this case $G$ is a thickening of a two-vertex trigraph and thus $G$ is a thickening of a trigraph in $\cal I$.

Suppose none of the last two cases occurs. Then $H_0$ is a block in which every edge is on a cycle of length $\ge2$. Let $s_0$ be the restriction of $s$ on $H_0$. Then $(H_0, s_0)$ is either in $\mathcal F_1\cup \mathcal F_2$ or evenly signed. First we assume $(H_0, s_0)$ is evenly signed. Then $(H,s)$ is also evenly signed. Moreover, $S_e$ is a thickening of a spring for every edge in $E_0=\{e\in E(H_0): s(e)=0\}$, and $S_e$ is a spot for every other edge of $H$. Let $S_e'$ be a spring for each $e\in E_0$ and let $S_e'=S_e$ for every other edge of $H$. Then $G$ is a thickening of $G'=\Omega(H,s,\{S_e'\})$. Now we only need to show that $G'\in \cal L$. Let $H'$ be obtained from $H$ by subdividing each edge in $E_0$ exactly once. Then $H'$ is bipartite. It follows from the construction of $\Omega$ that adjacent pairs of $G'$ are exactly adjacent pairs of the line graph $L(H')$. In addition, all semiadjacent pairs of $G'$ come from a spring, and thus no such pair is contained in a triangle. Therefore, $G'$ belongs to $\cal L$, as required.

It remains to consider the case $(H_0, s_0)\in \mathcal F_1\cup \mathcal F_2$. If $(H_0, s_0)\in \mathcal F_1$, then $H$ is obtained from $K_4$ by adding parallel edges and adding pendent edges to distinct vertices. Moreover, every $S_e$ is a spot. It follows that $G$ is an ordinary graph (meaning that $G$ has no semiadjacent pairs) and this graph is exactly $L(H)$. Now it is clear that $G$ is a thickening of $L(si(H))$, which belongs to $\mathcal J_1$. So we assume $(H_0, s_0)\in \mathcal F_2$. Let $V(H_0)=\{x_1, x_2, y_1, ..., y_m\}$ ($m\ge1$) such that $x_i$ ($i=1,2$) is adjacent to all other vertices. Like before, we assume that $H_0$ has no parallel edges, except for two possible edges $e_0,e_1$ between $x_1,x_2$, and such that $s(e_0)=0$ and $s(e_1)=1$. We also assume that $S_{e_0}$ is a spring, if $e_0$ is present. Suppose $H$ is obtained by adding pendent edges to $y_1,...,y_n$ ($n\ge0$) and to $k$ of $x_1,x_2$ ($0\le k\le 2$). If $e_0$ is present, then $G$ is a thickening of $J_2(n)$, where $\theta(uv)=0$. So assume that $e_0$ is not in $H$, and thus $G=L(H)$. For $i=1,2$, let $Q_i$ be the clique of $G$ formed by edges of $H$ incident with $x_i$. Let $Q_i'=Q_i-\{x_1x_2, x_iy_1,...,x_iy_n\}$. If $Q_1'\ne\emptyset$ is neither complete nor anticomplete to $Q_2'\ne\emptyset$, then again $G$ is a thickening of $J_2(n)$ with $\theta(uv)=0$. In the remainder cases (which are: some $Q_i'$ is empty, or  $Q_1'\ne\emptyset$ is complete or anticomplete to $Q_2'\ne\emptyset$), if $n=0$ then $G$ is a thickening of $K_3$, and if $n\ge 1$ then $G=J_2(n)-X$ for some $X\subseteq \{u,v\}$.  \qed

\section{Claw-free box-perfect graphs}

In this section we prove the following.

\begin{theorem} \label{thm:claw}
A claw-free perfect graph is box-perfect if and only if it is $S_3$-free.
\end{theorem}

We divide the proof into several lemmas. 
Let $G$ be a trigraph. 
We call $G$ a {\it sun} if $G\in tri(S_3)$. 
We call $G$ an {\it incomparability} trigraph if $G^{\ge0}$ is an incomparability graph. 
We call $G$ {\it elementary} if it is a thickening of a trigraph in $\cal L$. We remark that when an elementary trigraph has no semiadjacent pairs then they are exactly {\it elementary graphs} discussed in \cite{maffray}.

\begin{lemma} \label{lem:clawsun}
Let $G$ be a connected Berge trigraph. If $G$ is \{claw, sun\}-free then $G$ is obtained by simplicial summing incomparability trigraphs and elementary trigraphs. 
\end{lemma}

\vspace{-2mm}
\noindent{\bf Proof.} Since $G$ is connected, Berge, and claw-free, by Theorem \ref{thm:sum}, $G$ is obtained by simplicial summing thickenings of trigraphs in $\cal C\cup L\cup I\cup J$. Therefore, we may assume that $G$ is a thickening of a trigraph $G_0\in \cal C\cup L\cup I\cup J$. If $G_0\in \cal L$ then $G$ is elementary and we are done. If $G_0\in \cal C$ then $G_0|\{a_1,a_2,a_3,x_1,x_2,x_3\}$ (here we are using the notation in the definition of $\cal C$) is a sun and thus $G$ contains a sun, which is impossible. So we assume that $G_0\in \cal I\cup J$. In the following we prove that $G$ is an incomparability trigraph.

Suppose $G_0\in \cal I$. Then vertices of $G_0$ can be ordered as $v_1,...,v_n$ such that if $i<j<k$ and $\theta_0(v_iv_k)\ge0$ then $\theta_0(v_iv_j) = \theta_0(v_jv_k) =1$. Using the notation in the definition of thickening, we let $X_{v_i}=\{x_{i,j}:j=1,...,n_i\}$ ($1\le i\le n$). Now we define a binary relation $\prec$ on $V(G)$ such that $x_{i_1,j_1}\prec x_{i_2,j_2}$ if $\theta(x_{i_1,j_1} x_{i_2,j_2})=-1$ and $(i_1,j_1)$ is lexicographically smaller than $(i_2,j_2)$. We claim that $\prec$ is transitive. Suppose $x_{i_1,j_1}\prec x_{i_2,j_2} \prec x_{i_3,j_3}$. Since each $X_{v_i}$ is a strong clique, we must have $i_1<i_2<i_3$. It follows that $\theta_0(v_{i_1} v_{i_2})\le0$ and $\theta_0(v_{i_2} v_{i_3})\le0$. As a result, $\theta_0(v_{i_1} v_{i_3})=-1$ and thus $\theta(x_{i_1,j_1} x_{i_3,j_3})=-1$, which proves our claim. This claim implies that the complement of $G^{\ge0}$ is the comparability graph of poset $(V(G),\prec)$, which proves that $G^{\ge0}$ is an incomparability graph and thus $G$ is an incomparability trigraph.

Now suppose $G_0\in \cal J$. We claim that $G_0$ is a thickening of a trigraph in $\cal I$. This claim clearly implies that $G$ is a thickening of trigraph in $\cal I$, and thus the last paragraph proves that $G$ is an incomparability trigraph. 

Before proving the claim we make an observation. It is clear that every cobipartite trigraph is a thickening of a trigraph that has exactly two vertices and that the two vertices are semiadjacent. Since this two-vertex trigraph is in $\cal I$, our claim holds if $G_0$ is cobipartite.

We first consider the case $G_0\in \mathcal J_1$. If $G_0$ has two or more cubic vertices then $G_0$ contains an induced $S_3$. So $G_0$ contains at most one cubic vertex and in this case $G_0$ is cobipartite. Next we assume $G_0\in\mathcal J_2$ and let $G_0=J_2(n)-X$ (see the definition of $\mathcal J_2$). Let $x_iy_iz_i$ ($1\le i\le n$) denote the vertical triangles of $J_2(n)$, where $y_i\in Q_1$ and $z_i\in Q_2$. If $n\ge3$ then $G_0|\{w,x_1,y_1,z_1,y_2,z_3\}$ is a sun. So we have $n\le 2$. Now it is straightforward to verify that either $G_0$ is cobipartite, or $G_0$ contains a sun (found in a similar way), or $G_0=J_2(2)-\{u,v\}$. In the last case, $G_0$ is a thickening of the trigraph $G^*=(\{t_1,t_2,t_3,t_4,t_5\}, \theta^*)$, where $\theta^*(t_it_{i+1})=1$ ($i=1,2,3,4$), $\theta^*(t_2t_4)=0$, and $\theta^*(t_it_j)=-1$ for all other pairs. This completes the proof of our claim and also completes the proof of the lemma. \qed

Although simplicial sum was defined for trigraphs, this operation can be naturally inherited by ordinary graphs. Moreover, we have the following. 
 
\begin{lemma} \label{lem:simsum}
The simplicial sum of two ESP graphs is ESP.
\end{lemma}

\vspace{-2mm}
\noindent{\bf Proof.} 
Let $G$ be the simplicial sum of $G_1$ and $G_2$ over $x_1$ and $x_2$, where $G_i$ is an ESP graph with a simplicial vertex $x_i$ for $i=1, 2$.
Let $\Lambda$ be a set of maximal cliques of $G$. Note that
$N_{G_1}(x_1)\cup N_{G_2}(x_2)$ is the only maximal clique of $G$ that contains
edges between $N_{G_1}(x_1)$ and $N_{G_2}(x_2)$. For $i=1,2$, let $\Lambda_i$ consist of members of $\Lambda$ that
are cliques in $G_i$. Since $G_i$ is ESP, $\Lambda_i$ has an equitable subpartition
$(\Lambda_{i1}, \Lambda_{i2})$. If $\Lambda$ does not contain the clique $N_{G_1}(x_1)\cup N_{G_2}(x_2)$, then $(\Lambda_{11} \cup \Lambda_{21}, \Lambda_{12} \cup \Lambda_{22})$ is clearly an equitable subpartition of $\Lambda$.
Now assume that $\Lambda$ contains $N_{G_1}(x_1)\cup N_{G_2}(x_2)$. For $i=1,2$, let $\Lambda_i'=\Lambda_i\cup \{\{x_i\}\cup N_{G_i}(x_i)\}$. Note that $d_{\Lambda_i'}(x_i)=1$. 
Since $G_i$ is ESP, $\Lambda_i'$ has an equitable subpartition $(\Lambda_{i1}', \Lambda_{i2}')$. 
Without loss of generality, suppose $d_{\Lambda_{i1}'}(x_i)=1$ and $d_{\Lambda_{i2}'}(x_i)=0$ for $i=1, 2$.
Let $\Lambda '$ be obtained from $\Lambda_{11}' \cup \Lambda_{21}'$
by replacing the two cliques containing $x_1$ or $x_2$ by $N_{G_1}(x_1)\cup N_{G_2}(x_2)$; and set $\Lambda ''=\Lambda_{12}' \cup \Lambda_{22}'$.
Then $(\Lambda', \Lambda'')$ is an equitable subpartition of $\Lambda$. \qed

\begin{lemma} \label{lem:bipartite}
Let $\Lambda$ be a set of cliques of a graph $G$, for which $V(G)$ is partitioned into two cliques $X,Y$. Then $G$ has a multiset $\Lambda'$ of cliques such that \\ 
\indent (i) $|\Lambda'|=|\Lambda|$ and $d_{\Lambda'}(v)= d_\Lambda (v)$, for all $v\in V(G)$; \\ 
\indent (ii) members of $\Lambda'$ can be enumerated as $Q_1,...,Q_{|\Lambda|}$ such that every $v\in X$ appears in the first \\ \makebox[33pt]{} $d_{\Lambda}(v)$ terms and every $v\in Y$ appears in the last $d_{\Lambda}(v)$ terms.
\end{lemma}

\vspace{-2mm}
\noindent{\bf Proof.} For each $i=1,..., |\Lambda|$, let $X_i=\{x\in X: i\le d_{\Lambda}(x)\}$ and $Y_i=\{y\in Y: i\ge |\Lambda| - d_{\Lambda} (y) +1\}$. Then for every $i$, $Q_i=X_i\cup Y_i$ is a clique since $d_{\Lambda} (x) + d_{\Lambda}(y) \le |\Lambda|$ holds for all non-adjacent $x\in X$ and $y\in Y$. Now it is clear that $\Lambda'=\{Q_1,...,Q_{|\Lambda|}\}$ satisfies the requirements. 
\qed

\begin{lemma} \label{lem:elementary}
Elementary graphs are ESP.
\end{lemma}

\vspace{-2mm}
\noindent {\bf Proof.} Let elementary graph $H$ be obtained by thickening a trigraph $G$, where $G^{\ge0}$ is the line graph of a bipartite multigraph $B$ and such that semiadjacent pairs of $G$ are not contained in any triangle. Let $(Z_1, Z_2)$ be a partition of $V(B)$ into two stable sets. Let $u_1v_1, ..., u_nv_n$ be the semiadjacent pairs of $G$. Let $(X_v:v\in V(G))$ be the partition of $V(H)$ over which $G$ is thickened. For $i=1,...,n$, let $H_i = H[X_{u_i}\cup X_{v_i}]$. Since no semiadjacent pairs of $G$ are contained in a triangle, it is easy to see that for each maximal clique $C$ of $H$, either $C$ is a maximal clique of some $H_i$ or $C=\cup\{X_v:v\in Q\}$ for some maximal strong clique $Q$ of $G$. On the other hand, since $G^{\ge0}=L(B)$, for each maximal clique $Q$ of $G$ there exists a vertex $z$ of $B$ such that members of $Q$ are precisely edges of $B$ that are incident with $z$. We will say that $Q$ and $C=\cup\{X_v:v\in Q\}$ come from $z$. Note that, if $Q,Q'$ are maximal cliques of $G$ with $u_i\in Q-Q'$ and $v_i\in Q'-Q$ for some $i$, then $Q$ and $Q'$ come from vertices that both belong to $Z_1$ or both belong to $Z_2$.

Let $\Lambda$ be a set of maximal cliques of $H$. We need to show that $\Lambda$ admits an equitable subpartition. For $i=1,...,n$, let $\Lambda^{(i)}=\{C\in\Lambda: C\subseteq V(H_i)\}$. 
Let $\Lambda^{(0)} = \Lambda - \Lambda^{(1)} ... -\Lambda^{(n)}$. 
We assume by Lemma \ref{lem:bipartite} that members of each $\Lambda^{(i)}$ are enumerated as $C^{(i)}_1$, ...,  $C^{(i)}_{n_i}$ such that every $x\in X_{u_i}$ appears in the first $d_{\Lambda^{(i)}}(x)$ terms and every $x\in X_{v_i}$ appears in the last $d_{\Lambda^{(i)}}(x)$ terms. In the following we define a partition $(\Lambda_1,\Lambda_2)$ of $\Lambda$. To verify that $(\Lambda_1,\Lambda_2)$ is an equitable subpartition of $\Lambda$ we only need to verify $\min\{d_{\Lambda_1}(x), d_{\Lambda_2}(x)\}\ge \lfloor d_{\Lambda}(x)/2\rfloor$, for all $x\in V(H)$. 

We first consider $\Lambda^{(0)}$. If $C\in\Lambda^{(0)}$ then $C$ comes from a vertex $z$ of $B$. In this case we put $C$ into $\Lambda_i$ if $z\in Z_i$ ($i=1,2$). Since each $v\in V(G)$ is contained in at most two maximal cliques, we deduce that $d_{\Lambda}(x)\le 2$ for all $x\in V(H) -V(H_1)-... V(H_n)$. For these $x$, our partition of $\Lambda^{(0)}$ guarantees that $\min\{d_{\Lambda_1}(x), d_{\Lambda_2}(x)\}\ge \lfloor d_{\Lambda}(x)/2\rfloor$.

For cliques in each $\Lambda^{(i)}$ ($i=1,2,...,n$) we consider three cases. If none of $X_{u_i},X_{v_i}$ is contained in any clique of $\Lambda^{(0)}$, then $d_{\Lambda}(x)=d_{\Lambda^{(i)}}(x)$ for all $x\in V(H_i)$. Moreover, for each $x\in V(H_i)$, since cliques containing $x$ appear consecutively in the sequence $C^{(i)}_1$, ...., $C^{(i)}_{n_i}$, putting $C^{(i)}_j$ into $\Lambda_1$ for all odd $j$ and putting  $C^{(i)}_j$ into $\Lambda_2$ for all even $j$ lead to $\min\{d_{\Lambda_1}(x), d_{\Lambda_2}(x)\}\ge \lfloor d_{\Lambda}(x)/2\rfloor$.

If exactly one of $X_{u_i},X_{v_i}$ is contained in a clique of $\Lambda^{(0)}$, we assume by symmetry that $X_{u_i}$ is contained in a clique $C\in \Lambda^{(0)}$. We also assume without loss of generality that $C$ has been placed into $\Lambda_2$. For each $x\in V(H_i)$, since cliques containing $x$ appear consecutively in the sequence $X_{u_i},C^{(i)}_1$, ...., $C^{(i)}_{n_i}$, putting $C^{(i)}_j$ into $\Lambda_1$ for all odd $j$ and putting  $C^{(i)}_j$ into $\Lambda_2$ for all even $j$ lead to $\min\{d_{\Lambda_1}(x), d_{\Lambda_2}(x)\}\ge \lfloor d_{\Lambda}(x)/2\rfloor$.

If both $X_{u_i}, X_{v_i}$ are contained in cliques, say $C,D$, of $\Lambda^{(0)}$, by discussion in the first paragraph of this proof, we assume that both $C,D$ have been placed into $\Lambda_2$. 
%Note that $d_{\Lambda}(x) = 1 + d_{\Lambda^{(i)}}(x)$, for all $x\in V(H_i)$. 
We consider two subcases. Suppose $n_i$ is odd. For each $x\in V(H_i)$, since cliques containing $x$ appear consecutively in the sequence $X_{u_i},C^{(i)}_1$, ...., $C^{(i)}_{n_i},X_{v_i}$, putting $C^{(i)}_j$ into $\Lambda_1$ for all odd $j$ and putting  $C^{(i)}_j$ into $\Lambda_2$ for all even $j$ lead to $\min\{d_{\Lambda_1}(x), d_{\Lambda_2}(x)\}\ge \lfloor d_{\Lambda}(x)/2\rfloor$. Now suppose $n_i$ is even. For each $x\in V(H_i)$, note that cliques containing $x$ appear consecutively in the sequence $C^{(i)}_1,X_{u_i},C^{(i)}_2$, ...., $C^{(i)}_{n_i},X_{v_i}$, unless $x\in C_1^{(i)}\cap X_{v_i}$. In this case we put $C^{(i)}_j$ into $\Lambda_2$ for all odd $j>1$ and we put the rest into $\Lambda_1$. For each $x\in V(H_i) - (C_1^{(i)}\cap X_{v_i})$, it is clear that $\min\{d_{\Lambda_1}(x), d_{\Lambda_2}(x)\}\ge \lfloor d_{\Lambda}(x)/2\rfloor$. For each $x\in C_1^{(i)}\cap X_{v_i}$, we have $d_{\Lambda}(x)=1+n_i$. Our partition yields $d_{\Lambda_2}(x) = n_i/2$, which also leads to  $\min\{d_{\Lambda_1}(x), d_{\Lambda_2}(x)\}\ge \lfloor d_{\Lambda}(x)/2\rfloor$. \qed

\bigskip
\noindent {\bf Proof of Theorem \ref{thm:claw}.} The forward implication is clear, so we only need to consider the backward implication. Let $G$ be perfect and $\{claw, S_3\}$-free. By Lemma \ref{lem:clawsun}, each component of $G$ is obtained by simplicial summing incomparability graphs and elementary graphs. By Theorem \ref{thm:incomp} and Lemma \ref{lem:elementary}, incomparability graphs and elementary graphs are ESP. Thus $G$ is ESP by Lemma \ref{lem:simsum}, which proves that $G$ is box-perfect by Theorem \ref{thm:esp}. \qed

\end{document}